\documentclass{article}
\usepackage{amsmath,amsthm,amsfonts}

\numberwithin{equation}{section}

\def\ca{{\cal A}}
\def\cb{{\cal B}}
\def\cc{{\cal C}}
\def\cd{{\cal D}}
\def\ce{{\cal E}}
\def\cf{{\cal F}}
\def\cg{{\cal G}}
\def\ch{{\cal H}}

\def\cj{{\cal J}}
\def\ck{{\cal K}}

\def\cam{{\cal M}}
\def\cn{{\cal N}}
\def\co{{\cal O}}

\def\cu{{\cal U}}

\def\bn{{\mathbb N}}
\def\br{{\mathbb R}}
\def\bz{{\mathbb Z}}

\def\a{\alpha}
\def\b{\beta}
\def\g{\gamma}        \def\G{\Gamma}
\def\d{\delta}        \def\D{\Delta}
\def\eps{\varepsilon}

\def\th{\vartheta}

\def\l{\lambda}       
\def\m{\mu}
\def\n{\nu}

\def\p{\pi}
\def\r{\rho}

\def\t{\tau}

\def\f{\varphi}
\def\c{\chi}
\def\o{\omega}        \def\O{\Omega}

\def\imply{\Rightarrow}
\def\coimply{\Leftarrow}

\def\ov{\overline}
\def\e#1{{\rm e}^{#1}}

\def\limr{\lim_{r\to\infty}}
\def\limR{\lim_{R\to\infty}}
\def\lsup{\limsup_{R\to\infty}}
\def\linf{\liminf_{R\to\infty}}

\def\inj{{\rm inj}}
\def\supp{{\rm supp}\ }
\def\del{\bigtriangledown}

\def\ad#1{d_\infty(#1)}
\def\md#1{d_0(#1)}

\def\dE{{\partial E}}

\def\dec{\searrow}
\def\aff{\hat\in}
\def\nat{\natural}

\def\itm#1{\item[($#1$)]}

\newtheorem{Thm}{Theorem}[section]
\newtheorem{Cor}[Thm]{Corollary}
\newtheorem{Prop}[Thm]{Proposition}
\newtheorem{Lemma}[Thm]{Lemma}
\theoremstyle{definition}
\newtheorem{Dfn}[Thm]{Definition}
\newtheorem{exmp}[Thm]{Example}

\theoremstyle{remark}
\newtheorem{rem}[Thm]{Remark} 
\newtheorem{ack}{Acknowledgement} 
\title{\huge  Asymptotic dimension and\\
 Novikov-Shubin invariants\\
 for open manifolds}
\author{
D. Guido, T. Isola\\
Dipartimento di Matematica,\\ Universit\`a di Roma ``Tor
Vergata'',\\ I--00133 Roma, Italy.}

\date{December 23, 1996}
\begin{document}
\maketitle
\markboth{Asymptotic dimension package}
{Asymptotic dimension and Novikov-Shubin invariants}
\renewcommand{\sectionmark}[1]{}
\bigskip

 \begin{abstract}
A trace on the C$^*$-algebra $\ca$ of quasi-local operators on an open manifold
is described, based on the results in \cite{RoeOpen}. It allows a description
{\it \`a la} Novikov-Shubin \cite{NS2} of the low frequency behavior of the 
Laplace-Beltrami operator.
The 0-th Novikov-Shubin invariant defined in terms of such a trace is proved to
coincide with a metric invariant, which we call asymptotic dimension, thus
giving a large scale ``Weyl asymptotics'' relation.
Moreover, in analogy with the Connes-Wodzicki result \cite{CoCMP,Co,Wo}, the asymptotic
dimension $d$ measures the singular traceability (at 0) of the Laplace-Beltrami
operator, namely we may construct a (type II$_1$) singular trace which is
finite on the $^*$-bimodule over $\ca$ generated by $\Delta^{-d/2}$.
\end{abstract}

\newpage

 \setcounter{section}{-1}
 \section{Introduction.}\label{sec:intro}

 The inspiration of this paper came from the idea of Connes' 
 \cite{Co} of defining the dimension of a noncommutative compact manifold in
terms of the Weyl asymptotics, namely as the inverse of the
order of growth of the eigenvalues of differential operators of order one
(the Dirac operator for example).
Moreover Connes observed that a noncommutative measure (trace) may be attached
to such noncommutative dimension via the Dixmier trace, setting
$\t(a)={\rm{tr}}_\omega(a|D|^{-d})$, where $a$ is a ``function'' on the
noncommutative manifold, $D$ is the Dirac operator, $d$ is the noncommutative
dimension and ${\rm{tr}}_\omega$ is the (logarithmic) Dixmier trace.
According to the identification of the Dixmier trace with the Wodzicki residue,
such trace gives back the ordinary integration in the case of commutative
Riemannian manifolds.

In this paper we present a large scale analogue of these results for the case of
commutative noncompact manifolds.

We exhibit a large scale Weyl asymptotics, i.e. a correspondence between some
asymptotic dimension and the low frequency behavior of the Laplace-Beltrami
operator. Then we show that such asymptotic dimension carries a
(noncommutative) integration in terms of a type II$_1$ singular trace, which is
the low frequency analogue of the Dixmier trace.

The definition of asymptotic dimension is given in the context of metric
dimension theory, as a suitable large scale analogue of the metric dimension of
Kolmogorov and Tihomirov \cite{KT}.
Concerning the low frequency behavior, in case of a noncompact manifold
arising as universal  covering of a compact one, there is a set of numbers,
the so-called  Novikov-Shubin invariants $\a_{p}$, which are a measure of the
low frequency  behavior of the Laplacian (on $p$-forms) on the covering.

Since Atiyah, in his seminal paper \cite{Atiyah}, introduced a trace 
$\t$, called the $\G$-trace,
which replaced the ordinary trace in the statement of the index theorem for
covering manifolds $\G\to M \to X$ and brought to a definition of the Betti numbers for
coverings as $\t(\chi_{\{0\}}(\Delta))$,
Novikov and Shubin conjectured in \cite{NS1} that the behavior of
$\t(\chi_{[0,\l]}(\Delta))$ when $\l\to0$ should contain interesting topological
information.

Indeed the efforts of Novikov-Shubin \cite{NS2}, Lott \cite{Lott} and
Gromov-Shubin \cite{GS} proved that Novikov-Shubin numbers are indeed invariant
under homotopies of the base manifold.
The relations between Novikov-Shubin invariants and the singular traceability
of some $\Gamma$-invariant pseudodifferential operators is the object of a
separate paper \cite{GI3}.

In the case of open manifolds, J. Roe proved an index theorem \cite{RoeOpen} in
which he had to replace the $\Gamma$-trace of Atiyah (which, at least for
amenable coverings, may be seen as an average on the discrete group $\Gamma$)
with an ``average on the exhaustion'' trace.

We show that for open manifolds with bounded geometry and regular polynomial
growth the replacement of the Atiyah trace with the (suitably regularized)
Roe trace allows us to define the $0$-th Novikov-Shubin invariant for open
manifolds.

Then, the large scale Weyl asymptotics takes the form of the coincidence of the
asymptotic dimension with the $0$-th Novikov-Shubin invariant, which we prove
assuming the isoperimetric inequality of Grigor'yan \cite{Grigoryan94}.  Such a
relation shows in particular that the $0$-th Novikov-Shubin invariant, being a
metric object, is independent of  all the limiting procedures involved in its
definition.

The construction of the asymptotic (noncommutative) measure instead, depends 
on the singular traceability (at $0$) of $\Delta^{-\alpha/2}$, when the $0$-th
Novikov-Shubin invariant is a finite $\a\not=0$.

 For this we need more general singular traces, as singular traces on type I factors, as studied by Dixmier 
\cite{Dixmier}, Varga \cite{Varga} and  Albeverio-Guido-Ponosov-Scarlatti
\cite{AGPS}, only apply to compact  operators, like the negative powers of $\D$
on a compact manifold, but, in the case of non-type-I algebras, it has been
shown in \cite{GI1} the existence of singular traces which are finite on
suitable unbounded $\t$-compact operators, like negative powers of the
Laplacian on a noncompact manifold.

The theory of singular traces on C$^*$-algebras developed in \cite{GI2} may
then be used to construct a type II$_1$ singular trace on the unbounded
operators affiliated to a natural C$^{*}$-algebra of operators on the 
manifold.

 This paper is organized as follows.  \par 
 In section 1, a natural extension of the notion of 
 Kolmogorov-Tihomirov dimension \cite{KT} to nonnecessarily totally bounded 
 metric spaces is used, as a kind of analogy, to introduce an
 asymptotic dimension for metric spaces. It is proved that this dimension 
is invariant under rough isometries.  \par 
 In section 2, after some preliminaries on open manifolds of bounded 
 geometry and using recent estimates for the heat kernel by 
 Coulhon-Grigor'yan \cite{CG96}, we provide a relation between the asymptotic dimension and 
 the long time behavior of the heat kernel of the manifold 
 (Corollary \ref{2.1.7}), and establish a connection with 
 N.~Th.~Varopoulos' notion of asymptotic dimension for semigroups of 
 operators \cite{VSC}, as applied to the heat semigroup (Corollary \ref{2.2.6}); 
 finally we compare our definition with an analogous one given 
 recently by E.~B.~Davies \cite{Davies} for cylindrical ends.  \par 
 In section 3 we introduce the C$^{*}$-algebra of almost local 
 operators on a manifold of bounded geometry, as the norm closure of 
 the finite propagation speed operators, and show that $\cc_{0}$ 
 functional calculi of the Laplace-Beltrami operator are almost local 
 (Corollary \ref{3.1.6}); then, regularizing a previous construction 
 by J.~Roe \cite{RoeOpen}, we exhibit a weight on $\cb(L^{2}(M))$, 
 when $M$ is a manifold of bounded geometry and regular polynomial 
 growth (Proposition \ref{3.2.4}), which becomes 
 a semifinite, lower semicontinuous trace on the C$^{*}$-algebra of 
 almost local operators after two successive procedures of 
 regularization have been performed (note that these two procedures 
 are described abstractly) still retaining, though, the same value of 
 the original weight on the heat semigroup (Corollary \ref{3.4.3}).  \par 
 In the last section, after a brief exposition of the theory of 
 singular traces on C$^{*}$-algebras, which is the subject of a 
 separate publication, we define the ($0$-th) Novikov-Shubin invariant $\a_0(M)$ 
 for an open manifold $M$ of bounded geometry and regular polynomial growth 
 (Definition \ref{4.2.3}) and show (Corollary \ref{4.2.4}) an 
 asymptotic analogue of Wodzicki-Connes result, namely that 
 $\D^{-\a_{0}(M)/2}$ is singularly traceable at $0$, which is a 
 statement on the asymptotic behavior of the ``small eigenvalues'' of 
 the Laplacian. Here we observe that while in Wodzicki-Connes result 
 only logarithmic divergences appear, because manifolds are locally 
 regular, in our context  
 different divergences appear, and we recover the 
 logarithmic one in case of ``asymptotic regularity'', for example if 
 a discrete group acts on the manifold, \cite{GI3}. Finally, under more 
 restrictive hypotheses, we show that the Novikov-Shubin invariant 
 coincides with the asymptotic dimension of the manifold (Theorem 
 \ref{4.3.2}). This may be seen as a generalization of a result by
Varopoulos that $\a_0(M)= growth(\G)$, because of the rough-isometry 
invariance of the
asymptotic dimension, and Proposition \ref{2.1.4}.

 \section{Asymptotic dimension.}\label{sec:first} 

The main purpose of this section is the introduction of an asymptotic dimension
for metric spaces. To our knowledge, the notion of asymptotic dimension  in the
general setting of metric dimension theory has not been studied,  even though
Davies \cite{Davies} proposed a definition in the case of cylindrical ends of a
Riemannian manifold.

We shall give a definition of asymptotic dimension for a general metric space,
based on the (local) Kolmogorov dimension \cite{KT} and state its main
properties. We compare our definition with Davies' and also with the
notion of dimension at infinity for semigroups \cite{VSC} in
Section~\ref{sec:semigroupformula}.

 \subsection{Kolmogorov-Tihomirov metric dimension}
\label{subsec:KTdim}

In this subsection we recall a definition of metric dimension due
to Kolmogorov and Tihomirov \cite{KT}.
Quoting from their paper, a dimension ``corresponds to the possibility of
characterizing the ``massiveness'' of sets in metric spaces by the help of the
order of growth of the number of elements of their most economical 
$\eps$-coverings,
as $\eps\to0$''. Set functions retaining the
general properties of a dimension (cf. Theorem~\ref{Thm:Kdim}) have been
studied by several authors. Our choice of the
Kolmogorov dimension is due to the fact that it is suitable for
the kind of generalization we need in this paper, namely it quite naturally
produces a definition of asymptotic dimension.

In the following, unless otherwise specified, $(X,\d)$ will denote a metric
space, $B_X(x,R)$ the open ball in $X$ with centre $x$ and radius $R$,
$n_r(\O)$ the least number of open balls of radius $r$ 
which cover $\O\subset X$, and $\n_r(\O)$ the largest number of disjoint open 
balls of radius $r$ centered in $\O$.

The following lemma is proved in \cite{KT}. Due to some notational difference,
we include a proof.

 \begin{Lemma}\label{1.1.1} $n_r(X) \geq \n_r(X)\geq n_{2r}(X)$. 
 \end{Lemma}
 \begin{proof} 
We have only to prove the second inequality when $\n_r$ is finite.
 Let us assume that $\{B(x_i,r)\}_{i=1}^{\n_r(X)}$ are disjoint balls centered
in $X$ and observe that, for any $y\in X$, 
$\d(y,\bigcup_{i=1}^{\n_r(X)} B(x_i,r))<r$, otherwise $B(y,r)$ would be disjoint
from $\bigcup_{i=1}^{\n_r(X)} B(x_i,r)$, contradicting the maximality of
$\n_r$.
So for all $y\in X$ there is $j$ s.t. $\d(y, B(x_j,r))<r$, that is $X\subset
\bigcup_{i=1}^{\n_r(X)} B(x_i,2r)$, which implies the thesis.
 \end{proof}\medskip

 Kolmogorov and Tihomirov \cite{KT} defined a dimension for totally bounded 
metric spaces $E$ as
\begin{equation}\label{TKdim}
 \md{E}:=\limsup_{r\to0} \frac{\log n_r(E)}{\log(1/r)}.
\end{equation}
Then we may give the following definition.

\begin{Dfn}\label{Dfn:Kdim}  Let $(X,\d)$ be a metric space. Then, denoting by
$\cb(X)$ the family of bounded subsets of $X$, the metric Kolmogorov dimension
of $X$ is
 $$
 \md{X}:=\sup_{B\in\cb(X)}\limsup_{r\to0} \frac{\log n_r(B)}{\log(1/r)}.
 $$
 \end{Dfn}

Then the following proposition trivially holds.

 \begin{Prop}\label{1.1.3} If $\{B_n\}$ is an exhaustion of $X$ by bounded 
subsets, namely $B_n$ is increasing and for any bounded $B$ there exists $n$
such that $B\subseteq B_n$, one has $\md{X}=\lim_n\md{B_n}$. In particular,
 \begin{equation}\label{Kdimformula}
\md{X}=\lim_{R\to\infty}\limsup_{r\to0} \frac{\log
n_r(B_X(x,R))}{\log(1/r)}
 \end{equation}
\end{Prop}

\begin{rem} If bounded subsets of $X$ are not totally bounded, we could
define $\md{X}$ as the supremum over totally bounded subsets. These two
definitions, which agree e.g. on proper spaces, may be different in general.
For example an orthonormal basis in an infinite dimensional Hilbert space has
infinite dimension according to Definition~\ref{Dfn:Kdim}, but has zero
dimension in the other case. A definition of metric dimension which coincides
with $d_0$ on bounded subsets of $\br^p$ has been given by Tricot \cite{Tricot}
 in terms of rarefaction indices.
\end{rem}

 Let us now show that this set function satisfies the basic
properties of a dimension \cite{Pontriagin,Tricot}.

 \begin{Thm}\label{Thm:Kdim} The set function $d_0$ is a dimension, namely it
satisfies
\item{$(i)$}
 If $X\subset Y$ then $\md{X}\leq
\md{Y}$.
\item{$(ii)$} If $X_1,X_2\subset X$ then
$\md{X_1\cup X_2} = \max \{ \md{X_1}, \md{X_2} \}$.
\item{$(iii)$} If $X$ and $Y$ are metric spaces, 
then $\md{X\times Y} \leq \md{X}+\md{Y}$.
\end{Thm}

 \begin{proof}
Property $(i)$ easily follows from formula~(\ref{Kdimformula}).
\\
Now we prove $(ii)$.  The inequality 
 $\md{X_1\cup X_2}$ $\geq$ $\max \{ \md{X_1},\md{X_2} \}$ follows from $(i)$. 
 For the converse inequality, let $x_{i}\in X_{i}$, and set 
 $\d:=\d(x_{1},x_{2})$,
$d_1=\md{X_1}$, $d_2=\md{X_2}$, with e.g. $d_1\geq d_2$. If $d_1=\infty$ the
property is trivial, so we may suppose $d_1\in\br$.
Then
  $$
 B_{X_1\cup X_2}(x_1,R)\subset B_{X_1}(x_1,R)\cup B_{X_2}(x_2,R+\d)
 $$
therefore
 \begin{equation}\label{ineq1}
 n_r( B_{X_1\cup X_2}(x_1,R) )
  \leq n_r( B_{X_1}(x_1,R) ) + n_r( B_{X_2}(x_2,R+\d) ).
 \end{equation}
 Now, $\forall R>0$, 
$$
\limsup_{r\to0} \frac {\log n_r(B_{X_i}(x_{i},R))}{\log (1/r)}\leq  d_i
$$ 
i.e. $\forall R,\eps>0$ there is $r_0=r_0(\eps,R)$ s.t. for all $0<r<r_0$,
 $n_r(B_{X_1}(x_{1},R)) \leq r^{-(d_1+\eps)}$, and 
 $n_r(B_{X_2}(x_{2},R+\d)) \leq r^{-(d_2+\eps)}$
 hence, by (\ref{ineq1}),
 $$
 n_r(B_{X_1\cup X_2}(x,R))
 \leq r^{-(d_1+\eps)}(1+r^{d_1-d_2}).
 $$
 Finally,
 $$
 \limR\limsup_{r\to0} \frac {\log n_r(B_{X_1\cup X_2}(x,R))} {\log (1/r)} 
 \leq  d_1+\eps,
 $$
 that is
 $$
 \md{X_1\cup X_2}\leq
 \max\{\md{X_1},\md{X_2}\}+\eps
 $$
and the thesis follows by the arbitrariness of $\eps$.
\\
The proof of part $(iii)$ is postponed.
 \end{proof}\medskip

Kolmogorov dimension is indeed quasi-isometry invariant, as next proposition
shows.
 \begin{Prop}\label{1.1.8} Let $X,Y$ be metric spaces, and
$f:X\to Y$ a surjective quasi-isometry, namely $f$ satisfies
 $$
 c_1 \d_X(x_1,x_2) \leq \d_Y(f(x_1),f(x_2)) \leq c_2
\d_X(x_1,x_2).
 $$
 Then $\md{X}=\md{Y}$.
 \end{Prop}
 \begin{proof}
 By hypothesis we have $f(B_X(x,\r/c_2))\subset
B_Y(f(x),\r) \subset f(B_X(x,\r/c_1))$. So that, with
$y_j=f(x_j)$, $n:=n_r(B_Y(f(x),R))$,
 \begin{align*}
 f(B_X(x,R/c_2)) &\subset B_Y(f(x),R) \subset
 \bigcup_{j=1}^{n} B_Y(y_j,r) \\
 &\subset
 \bigcup_{j=1}^{n} f(B_X(x_j,r/c_1)) =
 f(\bigcup_{j=1}^{n} B_X(x_j,r/c_1))
 \end{align*}
 which implies $n_{r/c_1}(B_X(x,R/c_2))\leq
n_r(B_Y(f(x),R))$. 
\\
Since quasi-isometries are injective, we may repeat the same argument for
$f^{-1}$, and we get
$n_{c_2r}(B_Y(f(x),c_1R))
\leq n_r(B_X(x,R))$, so that
$n_{r/c_1}(B_X(x,R/c_2))$ $\leq$ 
$n_r(B_Y(f(x),R))$ $\leq$ $n_{r/c_2}(B_X(x,R/c_1))$.
 Finally 
\begin{align*}
 \limsup_{r\to0}\frac {\log n_{r/c_1}(B_X(x,R/c_2)) }{
 \log(c_1/r) - \log c_1}
 & \leq \limsup_{r\to0}\frac {\log n_r(B_Y(f(x),R)) }{
\log(1/r)} \\
 & \leq
\limsup_{r\to0}\frac {\log n_{r/c_2}(B_X(x,R/c_1)) }{
\log(c_2/r) - \log c_2}
 \end{align*} 
 and the thesis follows.
 \end{proof}\medskip

\noindent
{\it Proof of Theorem~\ref{Thm:Kdim}} (continued).
By the preceding Proposition, we may endow $X\times Y$
with any metric quasi-isometric to the product metric, i.e. 
$$
\d_{X\times Y}((x_1,y_1), (x_2,y_2))=\max\{ \d_X(x_1,x_2),\d_Y(y_1,y_2)\}.
$$ 
Then the thesis follows easily by 
$n_r(B_{X\times Y}((x,y),R))\leq n_r(B_X(x,R))\ n_r(B_Y(y,R))$. 
\qed\medskip

 \begin{rem}\label{1.1.9} Kolmogorov and Tihomirov assign a metric dimension to
a totally bounded metric space $E$ when $\exists\lim_{r\to\infty}$ in
equation (\ref{TKdim}), and consider upper and lower metric dimensions in the
general case. We observe that if the $\liminf$ is considered, the classical
dimensional inequality \cite{Pontriagin} stated in Theorem~\ref{Thm:Kdim} (iii) is
replaced by $\md{X\times Y} \geq \md{X}+\md{Y}$.
 \end{rem}

 \subsection{Asymptotic dimension}\label{subsec:asympdim}

 The function introduced in the previous subsection can be
used to study local properties of metric spaces. In
this paper we are mainly interested in the investigation of
the large scale behavior of these spaces, so we need a
different tool. Looking at equation~\ref{Kdimformula}, it is natural to set
the following

 \begin{Dfn}\label{1.2.1}  Let $(X,\d)$ be a  metric space. We call
 $$
\ad{X}:=\limr\lsup \frac{\log n_r(B_X(x,R)) }{ \log R},
 $$
 the {\it asymptotic dimension} of $X$. 
 \end{Dfn}

Let us remark that, as $n_r(B_X(x,R))$ is nonincreasing in
$r$, the function
 $$
r\mapsto\lsup\frac{\log n_r(B_X(x,R)) }{ \log R}
 $$ 
is nonincreasing too, so the $\lim_{r\to0}$ exists.

 \begin{Prop}\label{1.2.3} $\ad{X}$ does not depend
on $x$.
 \end{Prop} 
 \begin{proof}
 Let $x,y\in X$, and set $\d:=\d(x,y)$, so that $B(x,R)\subset
B(y,R+\d)
\subset B(x,R+2\d)$. This implies, 
\begin{align*}
 \frac{\log n_r(B(x,R)) }{ \log R} 
 &\leq
\frac{\log n_r(B(y,R+\d)) }{ \log (R+\d)}\
\frac{\log (R+\d) }{ \log R} \\
 &\leq
 \frac{\log n_r(B(x,R+2\d)) }{ \log (R+2\d)}\
 \frac{\log (R+2\d) }{ \log R}
 \end{align*}
 so that, taking $\lsup$ and then $\limr$ we
get the thesis.
 \end{proof}\medskip

 \begin{Lemma}\label{1.2.4}
  $$\ad{X}=\limr\lsup \frac{\log \n_r(B_X(x,R)) }{
\log R}$$ 
 \end{Lemma}
 \begin{proof} Follows easily from lemma \ref{1.1.1}.
 \end{proof}\medskip

 \begin{Thm}\label{Thm:adim} The set function $d_\infty$ is a dimension, 
namely it satisfies
\item{$(i)$}
 If $X\subset Y$ then $\ad{X}\leq
\ad{Y}$.
\item{$(ii)$} If $X_1,X_2\subset X$ then
$\ad{X_1\cup X_2} = \max \{ \ad{X_1}, \ad{X_2} \}$.
\item{$(iii)$} If $X$ and $Y$ are metric spaces, 
then $\ad{X\times Y} \leq \ad{X}+\ad{Y}$.
 \end{Thm}
 \begin{proof}
$(i)$ Let $x\in X$, then $B_X(x,R)\subset B_Y(x,R)$ and the
claim follows easily.
\\
$(ii)$  By part $(i)$, we get 
 $\ad{X_1\cup X_2} \geq \max \{ \ad{X_1},\ad{X_2} \}$. 
 Let us prove the converse inequality. 
\\
 Let $x_i\in X_i$, $i=1,2$, and set $\d=\d(x_1,x_2)$, $a=\ad{X_1}$,
$b=\ad{X_2}$, with e.g. $a\leq b$. 
 Then, $\forall \eps,r>0$  $\exists R_0=R_0(\eps,r)$ s.t. $\forall R>R_0$
 \begin{align*}
 n_r(B_{X_1}(x_1,R))  & \leq R^{a+\eps} \\
 n_r(B_{X_2}(x_2,R+\d)) & \leq R^{b+\eps},
 \end{align*}
 hence, by inequality (\ref{ineq1}),
 \begin{align*}
 n_r(B_{X_1\cup X_2}(x_1,R))
 &\leq R^{a+\eps}+R^{b+\eps}\\
 &= R^{b+\eps}(1+R^{a-b}).\\
 \end{align*}
 Finally,
 $$
 \frac {\log n_r(B_{X_1\cup X_2}(x_1,R))} {\log R}
 \leq b+\eps+ \frac {\log(1+R^{a-b})}{\log R}.
 $$
 Taking the $\lsup$ and then the $\limr$ we get
 $$
 \ad{X_1\cup X_2}\leq
 \max\{\ad{X_1},\ad{X_2}\}+\eps
 $$
and the thesis follows by the arbitrariness of $\eps$.
\\
The proof of part $(iii)$ is analogous to that of part $(iii)$ in
Theorem~\ref{Thm:Kdim}, where we may use Proposition~\ref{1.2.10}
because quasi-isometries are rough isometries. 
\end{proof}\medskip

\begin{rem}\label{1.2.7} In part $(ii)$ of the previous theorem we considered 
$X_1$ and $X_2$ as metric subspaces of $X$. If $X$ is a Riemannian manifold and
we endow the submanifolds $X_1$, $X_2$ with their geodesic metrics this
property does not hold in general. A simple example is the following. Let $f(t)
:= (t\cos t,t\sin t)$, $g(t) := (-t\cos t, -t\sin t)$, $t\geq0$ planar curves,
and set $X,\ Y$ for the closure in $\br^2$ of the two connected components of
$\br^2\setminus (G_f\cup G_g)$, where $G_f,\ G_g$ are the graphs of $f,\ g$,
and endow $X,\ Y$ with the geodesic metric. Then $X$ and $Y$ are
roughly-isometric to $[0,\infty)$ (see below) so that $\ad{X}=\ad{Y}=1$, while
$\ad{X\cup Y}=2$. 
 \end{rem}

 \begin{rem}\label{1.2.15} As for the local case, the choice of the 
 $\limsup$
in Definition~\ref{1.2.1} is the only one compatible with the classical
dimensional inequality stated in Theorem~\ref{Thm:adim} $(iii)$. This will
motivate our choice of the $\limsup$ in formula (\ref{e:NS}) for the 0-th
Novikov-Shubin invariant. 
 \end{rem} 

 \begin{Dfn}\label{1.2.8} Let $X,Y$ be  metric spaces,
$f:X\to Y$ is said to be a rough isometry if there are $a\geq1$,
$b,\eps\geq0$ s.t.
 \itm{i} $a^{-1}\d_X(x_1,x_2)-b \leq \d_Y(f(x_1),f(x_2)) \leq
a \d_X(x_1,x_2)+b$, for all $x_1,x_2\in X$, 
 \itm{ii} $\bigcup_{x\in X} B_Y(f(x),\eps) = Y$ 
 \end{Dfn}

It is clear that the notion of rough isometry is weaker then the notion of
quasi isometry introduced in the preceding subsection and, since any compact
set is roughly isometric to a point, $d_0$ is not rough-isometry invariant. 
We shall show that the asymptotic dimension is indeed invariant under rough 
isometries. 

 \begin{Lemma}\label{1.2.9} {\rm (\cite{Chavel}, Proposition 4.3)} If
$f:X\to Y$ is a rough isometry, there is a rough isometry
$f^-:Y\to X$, with constants $a,b^-,\eps^-$, s.t. 
 \itm{i} $\d_X(f^-\circ f(x),x)<c_X$, $x\in X$, 
 \itm{ii} $\d_Y(f\circ f^-(y),y)<c_Y$, $y\in Y$. 
 \end{Lemma}

 \begin{Prop}\label{1.2.10} Let $X,Y$ be  metric spaces, and
$f:X\to Y$ a  rough isometry. Then $\ad{X}=\ad{Y}$.
 \end{Prop}
 \begin{proof}
 Let $x_0\in X$, then for all $x\in B_X(x_0,r)$ we have
 $$
 \d_Y(f(x),f(x_0))\leq a \d_X(x,x_0)+b<ar+b
 $$ 
 so that
 $$
 f(B_X(x_0,r))\subset B_Y(f(x_0),ar+b).
 $$ 
 Then, with $n:=n_r(B_Y(f(x_0),aR+b))$, 
 $$
 f(B_X(x_0,R)) \subset 
 \bigcup_{j=1}^{n} B_Y(y_j,r),
 $$
 which implies 
 \begin{align*}
 f^-\circ f(B_X(x_0,R)) 
 & \subset \bigcup_{j=1}^{n} f^-(B_Y(y_j,r)) \\
 & \subset \bigcup_{j=1}^{n} B_X(f^-(y_j),ar+b^-).
 \end{align*}
 Let $x\in B_X(x_0,R)$, and $j$ be s.t. $f^-\circ f(x)\in
B_X(f^-(y_j),ar+b^-)$, then 
 $$
 \d_X(x,f^-(y_j))\leq \d_X(x,f^-\circ
f(x))+\d_X(f^-\circ f(x),f^-(y_j))<c_X+ar+b^-,
 $$ 
 so that  
 $$
 B_X(x_0,R)\subset
\bigcup_{j=1}^{n}
B_X(f^-(y_j),ar+b^-+c_X),
 $$
 which implies
$n_{ar+b^-+c_X}(B_X(x_0,R))\leq n_r(B_Y(f(x_0),aR+b))$. \\
 Finally 
 \begin{align*}
 \ad{X} 
 & = \limr\lsup\frac {\log n_r(B_X(x_0,R))}{\log R} \\
 & = \limr\lsup\frac {\log n_{ar+b^-+c_X}(B_X(x_0,R))}{\log R}\\
 & \leq \limr\lsup\frac {\log n_r(B_Y(f(x_0),aR+b))}{\log R} \\
 & = \limr\lsup\frac {\log n_r(B_Y(f(x_0),R))}{\log R} \\
 & = \ad{Y}
 \end{align*}
 and exchanging the roles of $X$ and $Y$ we get the thesis.
 \end{proof}\medskip

In what follows we show that when $X$ is equipped with a suitable measure, the
asymptotic dimension may be recovered in terms of the volume asymptotics for
balls of increasing radius, like the local dimension detects the volume
asymptotics for balls of infinitesimal radius.

 \begin{Dfn}\label{1.2.12} A Borel measure $\m$ on $(X,\d)$ is
said to be uniformly bounded if there are functions $\b_1,\b_2$,
s.t. $0<\b_1(r)\leq \m(B(x,r)) \leq
\b_2(r)$, for all $x\in X$, $r>0$. \\
 That is $\b_1(r):= \inf_{x\in X} \m(B(x,r)) >0$, and $\b_2(r) :=
\sup_{x\in X} \m(B(x,r)) <\infty$.
 \end{Dfn}

 \begin{Prop}\label{1.2.13} If $(X,\d)$ has a uniformly
bounded measure, then every ball in $X$ is
totally bounded (so that if $X$ is complete it is locally compact).
 \end{Prop}
 \begin{proof}
Indeed, if there is a ball $B=B(x,R)$ which is not totally
bounded, then there is $r>0$ s.t. every $r$-net in $B$ is
infinite, so $n_r(B)$ is infinite, and $\n_r(B)$ is infinite
too. So that $\b_2(R)\geq \m(B) \geq \sum_{i=1}^{\n_r(B)}
\m(B(x_i,r)) \geq \b_1(r)\n_r(B) = \infty$, which is absurd.
 \end{proof}\medskip

 \begin{Prop}\label{1.2.14} If $\m$ is a uniformly bounded Borel
measure on $X$ then 
 $$
 \ad{X}=\lsup\frac{\log\m(B(x,R))}{ \log R} .
 $$ 
 \end{Prop}
 \begin{proof}
 As $\bigcup_{i=1}^{\n_r(B(x,R))} B(x_i,r) \subset B(x,R+r)
\subset \bigcup_{j=1}^{n_r(B(x,R+r))} B(y_j,r)$, we get
 $$
 \b_2(r)n_r(B(x,R+r))\geq \m(B(x,R+r)) \geq
\b_1(r)\n_r(B(x,R)) \geq \b_1(r)n_{2r}(B(x,R)),
 $$
 by Lemma \ref{1.1.1}. So that 
 $$
 \b_1(r/2)\leq
\frac{\m(B(x,R+r/2))}{ n_r(B(x,R))},
\qquad
\frac{\m(B(x,R))}{ n_r(B(x,R))} \leq \b_2(r),
 $$
 and the thesis follows easily.
 \end{proof}\medskip

Let us conclude this subsection with some examples.

 \begin{exmp}\label{1.2.16} 
 \itm{i} $\br^n$ has asymptotic dimension $n$. 
 \itm{ii} Set $X:= \cup_{n\in\bz}\{(x,y)\in\br^2 :
\d((x,y),(n,0))<\frac{1}{4} \}$, endowed with the Euclidean
metric, then $\md{X}=2$, $\ad{X}=1$. 
\itm{iii} Set $X=\bz$ with the counting measure,
then $\md{X}=0$, and $\ad{X}=1$. 
\itm{iv} Let $X$ be the unit ball in an infinite dimensional Banach space. Then
$\md{X}=+\infty$ while $\ad{X}=0$.
 \end{exmp}

 \begin{exmp}\label{1.2.17} Set $X:=\{(x,y)\in\br^2: x\geq0, |y|\leq x^\a\}$,
endowed with the Euclidean metric, where $\a\in(0,1]$. Then
$\ad{X}=\a+1$.
 \end{exmp}
 \begin{proof}
 This metric space has a uniformly bounded Borel measure, the
Lebesgue area, so we can use Proposition \ref{1.2.14}.
 Set $x_0:=(0,0)$, and $B_R:= B_X(x_0,R)$. Then, if $R\geq
\sqrt{4^{1+1/\a}r^{2/\a}+4^{1+\a}r^2}$, $B_{R}\subset Q_{1}\cup 
Q_{2}$, where $Q_1:= \{ (x,y)\in\br^2 : -2r\leq x \leq (2r)^{1/\a}+2r,\ 
|y|\leq 4r \}$, 
and $Q_2:=\{ (x,y)\in\br^2 : (2r)^{1/\a} \leq x \leq R,
\ |y|\leq 2x^\a \}$. Now, if $x_R>0$ is s.t.
$x_R^2+x_R^{2\a}=R^2$, we get
 \begin{align*}
 area(B_R) & \geq \frac2{\a+1} x_R^{\a+1} \\
 area(Q_1) & = 4r(4r+(2r)^{1/\a}) \\
 area(Q_2) & = \frac4{\a+1} (R^{1+\a}-(2r)^{1+1/\a}),
 \end{align*}
 so that
 \begin{align*}
 \limR \frac {(\a+1)\log x_R}{\log R} 
 & \leq \linf \frac {\log area(B_R)}{\log R} \\
 & \leq \lsup \frac {\log area(B_R)}{\log R} \leq \a+1
 \end{align*}
 and, as $\limR \frac {\log x_R}{\log R} = \lim_{x\to\infty}
\frac {\log x}{\log\sqrt{x^2+x^{2\a}}}=1$, we get the thesis.
 \end{proof}\medskip

 \section{A semigroup formula for the asymptotic dimension of
an open manifold}\label{sec:semigroupformula}

 \subsection{Open manifolds of bounded 
 geometry}\label{subsec:openmanifold}
 
 In this subsection, after some preliminary results on open
manifolds of bounded geometry, we give a formula for the
asymptotic dimension in terms of the asymptotics of the
heat kernel. This opens the way for the abstract treatment
of the following subsection.  \par 
 Several definitions of bounded geometry for an open manifold (i.e. a
Riemannian, complete, noncompact manifold) are
usually considered. They all require some uniform bound (either
from above or from below) on some geometric objects, such as:
injectivity radius, sectional curvature, Ricci curvature, Riemann
curvature tensor etc. (For all unexplained notions see e.g. Chavel's
book \cite{Chavel}).
 \par 
 In this paper the following form is used, but see \cite{RoeOpen} and
references therein for a different approach.  \par 

 \begin{Dfn}\label{2.1.1}  Let $(M,g)$ be an
$n$-dimensional complete Riemannian manifold. We say that $M$ has
bounded geometry if it has positive injectivity radius, sectional
curvature bounded from above by some constant $c_1$, and Ricci
curvature bounded from below by $(n-1)c_2 g$. 
 \end{Dfn}
 
 \begin{Thm}\label{2.1.2} {\rm (\cite{Chavel}, p.119,123)} Let $M$ be a
complete Riemannian manifold of bounded geometry. Then there are
positive real functions $\b_1,\ \b_2$ s.t.
 \itm{i}
 $$
 0<\b_1(r)\leq vol(B(x,r)) \leq \b_2(r),
 $$
 for all $x\in X$, $r>0$, that is the volume form is a uniformly
bounded measure (cf. Definition \ref{1.2.12}), 
 \itm{ii} $\lim_{r\to0}\frac{\b_2(r)}{\b_1(r)} =1$. 
 \end{Thm}
 \begin{proof}
 $(i)$ We can assume $c_2<0<c_1$ without loss of generality. Then,
denoting with $V_\d(r)$ the volume of a ball of radius $r$ in a
manifold of constant sectional curvature equal to $\d$, we can set
 $\b_1(r) := V_{c_1}(r\wedge r_0)$, and $\b_2:= V_{c_2}(r)$,
where $r_0:= \min\{ \inj(M), \frac{\pi}{\sqrt{c_1}} \}$, and 
$\inj(M)$ is the injectivity radius of $M$. \\
 $(ii)$ 
 \begin{align*}
 \lim_{r\to0} \frac{\b_2(r)}{\b_1(r)}  
 & = \lim_{r\to0} \frac{V_{c_2}(r)}{V_{c_1}(r)} 
   = \lim_{r\to0} \frac{ \int_0^r S_{c_2}(t)^{n-1}dt }{
                    \int_0^r S_{c_1}(t)^{n-1}dt } \\
 & = \left( \lim_{r\to0} \frac{S_{c_2}(r)}{S_{c_1}(r)} \right)^{n-1} = 1
 \end{align*}
	 where (cfr. \cite{Chavel}, formulas (2.48), (3.24), (3.25)) $V_\d(r) =
\frac{n\sqrt{\pi} }{ \G(n/2+1)} \int_0^r S_\d(t)^{n-1}dt$, and
 $$
  S_\d(r) := 
  \begin{cases}
 \frac1{\sqrt{-\d}}\sinh(r\sqrt{-\d}) & \d<0 \\
        r& \d=0 \\
 \frac1{\sqrt{\d}}\sin(r\sqrt{\d}) &  \d>0. 
   \end{cases}
 $$
 \end{proof}\medskip 
 
Conditions under which the inequality in Theorem~\ref{Thm:adim} $(iii)$ becomes an
equality are often studied in the case of (local) dimension theory (cf.
\cite{Pontriagin,Salli}). The following proposition gives such a condition for
the asymptotic dimension.

 \begin{Prop}\label{2.1.3} Let $M,N$ be complete Riemannian
manifolds of bounded geometry, which admit asymptotic
dimension in a strong sense, that is
 $$
 \ad{M} \equiv \limr\limR \frac{\log n_r(B_M(o,R)) }{ \log R},
 $$
and analogously for $N$. Then 
 $$
 \ad{M\times N}=\ad{M}+\ad{N}.
 $$
 \end{Prop}
 \begin{proof}
 As $vol(B_{M\times N}((x,y),R))=vol(B_M(x,R))
vol(B_N(y,R))$, we get 
 \begin{align*}
\ad{M\times N} &= 
\limR \frac{\log vol (B_{M\times N}((x,y),R))}{ \log R} \\
 & = 
\limR \frac{\log vol (B_M(x,R))}{ \log R} + 
\limR \frac{\log vol (B_N(y,R))}{ \log R} \\ 
 & = \ad{M}+\ad{N}.
 \end{align*}
 \end{proof}\medskip

 As the asymptotic dimension is invariant under rough isometries,
it is natural to substitute the continuous space with a coarse
graining, which destroys the local structure, but preserves the
large scale structure. To state it more precisely, recall
(\cite{Chavel}, p. 194) that a discretization of a metric space $M$
is a graph $G$ determined by an $\eps$-separated subset $\cg$ of
$M$ for which there is a $R>0$ s.t. $M=\cup_{x\in\cg} B_M(x,R)$.
 The graph structure on $\cg$ is determined by one oriented edge
from any $x\in\cg$ to any $y\in\cg$, $y\neq x$, denoted $<x,y>$,
precisely when $\d_M(x,y)<2R$. Define the combinatorial metric on
$G$ by $\d_c(x,y):=\inf \{ \sum_{i=0}^n \d(x_i,x_{i+1}):
(x_0,\ldots,x_{n+1})\in Path_n(x,y),\ n\in\bn\}$, where
$Path_n(x,y) :=\{ (x_0,\ldots,x_{n+1}) : x_i\in\cg,\ x_0=x,\
x_{n+1}=y, <x_i,x_{i+1}>\in G\}$.  \par 

\begin{Prop}\label{2.1.4} {\rm (\cite{Chavel}, Theorem 4.9)} 
 Let $M$ be a complete Riemannian manifold with Ricci curvature
bounded from below. Then $M$ is roughly isometric to any of its
discretizations, endowed with the combinatorial metric. Therefore
$M$ has the same asymptotic dimension of any of its
discretizations.
 \end{Prop}

The previous result, together with the invariance of the asymptotic dimension
under rough isometries, shows that, when $M$ has a discrete group of isometries
$\Gamma$ with a compact quotient, the asymptotic dimension of the manifold
coincides with the asymptotic dimension of the group, hence with its growth
(cf. \cite{GI3}), hence, by the result of Varopoulos \cite{Varopoulos2}, 
it coincides with the 0-th Novikov-Shubin invariant.

 Let $M$ be a complete Riemannian manifold, and recall
(\cite{Daviesbook}, Chapter 5) that $\D$, the Laplace-Beltrami
operator on $M$, is essentially self-adjoint and positive on
$L^2(M)$, and the semigroup $\e{-t\D}$ has a strictly positive
$\cc^\infty$ kernel, $p_t(x,y)$, on $(0,\infty)\times M\times M$,
called the heat kernel. Recall the following results, which
will be useful in the sequel, and where we use $V(x,r):= vol(B(x,r))$, 
for simplicity.

 \begin{Prop}\label{2.1.5} Let $M$ be an $n$-dimensional complete Riemannian
manifold with Ricci curvature bounded below, and let $E$ be the 
infimum of the spectrum of $-\D$, then for all $\eps>0$, there are 
$c, c'>0$, s.t.
 \itm{i} 
 $$
 p_{t}(x,y)\leq 
 \begin{cases}
 c\  V(x,\sqrt{t})^{-1/2} V(y,\sqrt{t})^{-1/2} 
 \e{\frac{-\d(x,y)^{2}}{(4+\eps)t}} & 0<t\leq1 \\
 c\  V(x,1)^{-1/2} V(y,1)^{-1/2} \e{(\eps-E)t} 
 \e{\frac{-\d(x,y)^{2}}{(4+\eps)t}} & t\geq1
 \end{cases}
 $$
 \itm{ii}
  $$
 |\del_{x}p_{t}(x,y)|\leq 
 \begin{cases}
 c' t^{-n/2} (t^{-1}+\d(x,y))^{1/2}
 \e{\frac{-\d(x,y)^{2}}{(4+\eps)t}} & 0<t\leq1 \\
 c'  \e{(\eps-E)t} (t^{-1}+\d(x,y))^{1/2}
 \e{\frac{-\d(x,y)^{2}}{(4+\eps)t}} & t\geq1
 \end{cases}
 $$
 As a consequence $p_t$ is uniformly continuous on a neighborhood of the 
 diagonal of $M\times M$.
 \end{Prop}
 \begin{proof}
 $(i)$ is  (\cite{Davies88}, theorems 16, 17). \\
 Observe now that, with $r_{1}:= \min\{1,\inj{(M)}, 
 \frac\pi{\sqrt{c_{1}}}\}$, we have $V(x,t)\geq V(x,r_{1})\geq 
 V_{c_{1}}(r_{1})$ for any $t\geq1$. On the other hand, if 
 $t\in[r_{1},1]$, we have 
 $\frac{V(x,t)}{t^{n}}\geq V(x,r_{1})\geq V_{c_{1}}(r_{1})$, while, if 
 $t\in(0,r_{1})$, from $\frac{V(x,t)}{t^{n}}\geq 
 \frac{V_{c_{1}}(t)}{t^{n}}$ and $\lim_{t\to0} 
 \frac{V_{c_{1}}(t)}{t^{n}} = \frac{\sqrt{\pi}}{\G(n/2+1)}>0$ (use 
 the formulas in the proof of Theorem \ref{2.1.2}), there follows $a>0$ 
 s.t. $\frac{V(x,t)}{t^{n}}\geq a$ for any $t\in(0,r_{1})$. \\
 Putting all things together we get a simplified version of the 
 estimates $(i)$
  $$
 p_{t}(x,y)\leq 
 \begin{cases}
 c t^{-n/2} 
 \e{\frac{-\d(x,y)^{2}}{(4+\eps)t}} & 0<t\leq1 \\
 c \e{(\eps-E)t} 
 \e{\frac{-\d(x,y)^{2}}{(4+\eps)t}} & t\geq1\ .
 \end{cases}
 $$
 $(ii)$ follows from  (\cite{DaviesJOT}, theorem 6), using the 
 simplified estimates above. \\
 Finally, for any $\d_{0}<r_{1}$, $x\in M$, $y\in B(x,\d_{0})$, we 
 have $|p_{t}(x,y)-p_{t}(x,x)|\leq \sup|\del_{y}p_{t}(x,y)| \d(x,y)$, 
 and from $(ii)$ we get the uniform continuity.
 \end{proof}\medskip

 The following proposition shows the deep connection between the
heat kernel and the volume of balls.

 \begin{Thm}\label{2.1.6} {\rm (\cite{CG96}, Corollary 7.3)
(\cite{Grigoryan94}, Proposition 5.2)} \\
 Let $M$ be a complete Riemannian manifold, and set $\l_1(U)$ for
the first Dirichlet eigenvalue of $-\D$ in $U$. Then the
following are equivalent
 \itm{i} there are $\a,\ \b>0$ s.t. for all $x\in M$, $r>0$, and
all regions $U\subset B(x,r)$,  \\
 $$
 \l_1(U) \geq \frac{\a }{ r^2} \left( \frac{V(x,r)}{ vol(U)}
\right)^\b
 $$ 
 \itm{ii} there are $A,\ C,\ C'>0$ s.t. for all  $x\in M$,
$r>0$,
 \begin{align}\label{e:voldouble}
 V(x,2r) & \leq A V(x,r) \\
 \frac{C}{ V(x,\sqrt{r})} & \leq p_r(x,x) \leq 
 \frac{C'}{V(x,\sqrt{r})} . \notag 
 \end{align}
 \end{Thm}

 Following \cite{CG96} we call (\ref{e:voldouble}) the volume doubling
property.

 As a consequence of the above results we have the following

 \begin{Cor}\label{2.1.7} Let $M$ be a complete Riemannian manifold
of bounded geometry, and assume one of the equivalent properties of the
previous Theorem. \\ 
Then $\ad{M}=\limsup_{t\to\infty} \frac{-2\log p_t(x_0,x_0)}{\log t}$
 \end{Cor}
 \begin{proof} Follows from Theorem \ref{2.1.2}.
 \end{proof}\medskip

 For a different result under weaker hypotheses, see Corollary
\ref{2.2.6}.

 Before closing this subsection we observe that the volume
doubling property is a key notion for our work, as it is a weak
form of polynomial growth condition, and still guarantees the
finiteness of the asymptotic dimension (for manifolds of bounded
geometry).

 \begin{Prop}\label{2.1.8} Let $M$ be a complete Riemannian
manifold of bounded geometry, and suppose the volume doubling
property (\ref{e:voldouble}) holds. Then $\ad{M}\leq \log_2 A$.
 \end{Prop}
 \begin{proof}
 Let $R>1$, and $n\in\bn$ be s.t. $2^{n-1}<R \leq
2^n$. Then $V(x,R)\leq V(x,2^n)\leq A^n V(x,1)$, so that
 $$
 1 \leq \frac{V(x,R) }{ V(x,1)} \leq A^n \leq A R^{\log_2 A}.
 $$
 Therefore $\ad{M} = \limsup_{R\to\infty} \frac{\log V(x,R) }{ \log R} \leq \log_2
A$.
 \end{proof}\medskip

 \subsection{Asymptotic dimension of some semigroups of 
bounded operators}\label{subsec:asympsemi}

 Based on the notion of dimension at infinity due to
 Varopoulos, Saloff-Coste, Coulhon \cite{VSC}, see also \cite{Cou}, we define the
asymptotic dimension of a semigroup of bounded operators on a
measure space.

 \begin{Dfn}\label{2.2.1} Let $(X,\O,\m)$ be a measure space, and
$T_t : L^1(X,\O,\m)\to L^\infty(X,\O,\m)$ be a semigroup of
bounded operators. Then we set
 $$
 \ad{T} := \liminf_{t\to\infty}  \frac{-2\log \|T_t\|_{1\to\infty}
}{ \log t} .
 $$ 
 \end{Dfn}

 \begin{Thm}\label{2.2.2} {\rm (\cite{VSC}, Theorem II.4.3)}\\
 Let $(X,\O,\m)$ be a measure space, and suppose given 
$T_t\in \cb(L^1(X,\O,\m)\cap L^\infty(X,\O,\m))$, which, for any
$p\in[1,\infty]$, extends to a semigroup on $L^p$, of class
$C^0$ if $p<\infty$. Let $A$ be the generator, and suppose that
$T_t$ is equicontinuous on
$L^1$ and $L^\infty$, bounded  analytic on $L^2$, and
$\|T_1\|_{1\to\infty}<\infty$. Let $0<\a<\frac{n}{2}$. Then the
following are equivalent 
 \itm{i} $\|f\|_{2n/(n-2\a)} \leq C(\|A^{\a/2}f\|_2 +
\|A^{\a/2}f\|_{2n/(n-2\a)})$, $f\in \cd$ 
 \itm{ii} $\|T_1f\|_{2n/(n-2\a)} \leq C \|A^{\a/2}f\|_2$, $f\in
\cd$ 
 \itm{iii} $\|T_t\|_{1\to\infty} \leq Ct^{-n/2}$,
$t\in[1,\infty)$, \\
 where $\cd:= {\rm span\ } \{ \int_0^\infty \f(t)T_tf dt\ :\
\f\in\cc^\infty_0(0,\infty),\ f\in\ L^\infty(X,\O,\m),\
\m\{f\neq0\}<\infty \}$.
 \end{Thm}

 \begin{Prop}\label{2.2.3} Let $\{T_t\}$ be as in the previous
Theorem. Then the following are equivalent
 \itm{i} $\|T_t\|_{1\to\infty} \leq Ct^{-n/2}$, $t\geq1$
 \itm{ii} $\|T_t\|_{1\to\infty} \leq Ct^{-n/2}$, $t\geq t_0>1$. 
 \end{Prop}
 \begin{proof}
 $(ii)\imply(i)$ \\
 Let $t>1$ and observe that $\|T_t\|_{1\to\infty} =
\|T_1T_{t-1}\|_{1\to\infty} \leq
\|T_1\|_{1\to\infty}\|T_{t-1}\|_{1\to1} \leq
k\|T_1\|_{1\to\infty} =: M$, where $k:= \sup_{t>0}
\|T_t\|_{1\to1}<\infty$ because $T_t$ is equicontinuous on $L^1$.
 So that, with $C_0:= \max\{ C, Mt_0^{n/2} \}$, we get the thesis.
 \end{proof}\medskip

 \begin{Prop}\label{2.2.4} $\ad{T} = \sup \{n>0:
\|T_t\|_{1\to\infty} \leq Ct^{-n/2},\ t\geq1 \}$.
 \end{Prop}
 \begin{proof}
 Set $d$ for the supremum. Then for all $\eps>0$, there is $t_0>1$
s.t. $\|T_t\|_{1\to\infty}\leq t^{-(\ad{T}-\eps)/2}$, for all
$t\geq t_0$, and, by previous proposition, $\ad{T}-\eps\leq d$.
 Conversely $\|T_t\|_{1\to\infty}\leq t^{-(d-\eps)/2}$, for all
$t\geq1$ implies $d-\eps\leq\ad{T}$. 
 \end{proof}\medskip

  Using a recent result by Coulhon-Grigor'yan \cite{CG96} we can
show the relation between the asymptotic dimension of the heat
kernel semigroup and the asymptotic dimension of the
underlying manifold.

 \begin{Thm}\label{2.2.5} {\rm (\cite{CG96}, Theorem 2.7)} \\
Let $M$ be a complete Riemannian manifold. If there are
$x_0\in M$, $A>0$  s.t. 
$V(x_0,2t)\leq A V(x_0,t)$, for all $t>0$, then there is $c>0$ s.t.
 $$
 \sup_{x\in M} p_t(x,x) \geq \frac{c}{ V(x_0,\sqrt{t})}\ .
 $$
 \end{Thm} 

 \begin{Cor}\label{2.2.6} Let $M$ be a complete Riemannian manifold
of bounded geometry, and assume there are $x_0\in M$, $A>0$  s.t. 
$V(x_0,2t)\leq A V(x_0,t)$, for all $t>0$. 
 Then $\ad{\e{-t\D}}\leq \ad{M}$.
 \end{Cor}
 \begin{proof}
 Follows immediately from the previous results if we recall that
$\|\e{-t\D}\|_{1\to\infty}=\sup_{x\in M} p_t(x,x)$.
 \end{proof}\medskip

 \subsection{Asymptotic dimension of some cylindrical ends}
 \label{subsec:cylindrical}

 In this subsection we want to compare our work with a
recent work of E.B.~Davies'. In \cite{Davies} he defines the
asymptotic dimension of cylindrical ends of a Riemannian manifold $M$ as
follows. Let $E\subset M$ be homeomorphic to $(1,\infty)\times A$,
where $A$ is a compact Riemannian manifold. Set
$\dE:=\{1\}\times A$, $E_r:=\{ x\in E: \d(x,\dE) < r\}$, where $\d$
is the restriction of the metric in $M$. Then $E$ has asymptotic
dimension $D$ if there is a positive constant $c$ s.t. 
 \begin{equation}\label{e:davies}
  c^{-1} r^D \leq vol(E_r) \leq cr^D,
 \end{equation}
 for all $r\geq1$. He does not assume
bounded geometry for $E$. If one does, the two definitions
coincide as in the following

 \begin{Prop}\label{2.3.1} With the above notation, if the volume 
 form on $E$ is a uniformly bounded measure (as in Definition 
 \ref{1.2.12}), or in particular if
$E$ has bounded geometry (as in Definition \ref{2.1.1}), and there 
is $D$ as in {\rm (\ref{e:davies})}, then
$\ad{E}=D$. 
 \end{Prop}
 \begin{proof}
 Choose $o\in E$, and set
$\d:=\d(o,\dE)$, $\D:=diam(\dE)$. Then it is easy to
prove that $E_{R-\d-\D} \subset B_E(o,R) \subset E_{R+\d}$.
 \\
 Then $c^{-1}(R-\d-\D)^D \leq vol(B_E(o,R)) \leq c(R+\d)^D$,
and from \ref{1.2.14} the thesis follows. 
 \end{proof}\medskip

 Motivated by (\cite{Davies}, example 16), let us set the
following

 \begin{Dfn}\label{2.3.2} Let us say that $E$ is a standard end 
of local dimension $N$ if it is
homeomorphic to $(1,\infty)\times A$, endowed with the metric
$ds^2=dx^2+f(x)^2d\o^2$, and with the volume form
$dvol=f(x)^{N-1}dxd\o$, where $(A,\o)$ is an
$(N-1)$-dimensional compact Riemannian manifold, and $f$ is an
increasing smooth function. 
 \end{Dfn}

 \begin{Prop}\label{2.3.3} The volume form on a standard end $E$ is
a uniformly bound\-ed measure. Therefore, if $E$ satisfies {\rm (\ref{e:davies})},
we get $\ad{E}=D$.
 \end{Prop}
 \begin{proof}
 It is easy to show that, for $(x_0,p_0)\in E$, $r<x_0-1$, \\
 \begin{align*}
 [x_0-r/2,x_0+r/2]&\times B_A\left(p_0, \frac{r/2}{f(x_0+r/2)}\right) 
 \subset B_E((x_0,p_0),r) \\
 & \subset [x_0-r,x_0+r]\times B_A\left(p_0, \frac{r}{f(x_0-r)}\right)
 \end{align*}
 So that 
 \begin{align*}
 \int_{x_0-r/2}^{x_0+r/2}f(x)^{N-1}dx\ & V_A\left(p_0,
 \frac{r/2}{f(x_0+r/2)}\right) 
 \leq V_E((x_0,p_0),r) \\
 &\leq \int_{x_0-r}^{x_0+r}f(x)^{N-1}dx\ V_A\left(p_0,
  \frac{r}{f(x_0-r)}\right)
 \end{align*}
 which implies
 \begin{align*}
 rf(x_0-r/2)^{N-1}\ V_A\left(p_0,\frac{r/2}{f(x_0+r/2)}\right)
 &\leq V_E((x_0,p_0),r) \\
 &\leq 2rf(x_0+r)^{N-1}\ V_A\left(p_0,\frac{r}{f(x_0-r)}\right)
 \end{align*}
 As for $x_0\to\infty$, $V_A(p_0,\frac{r}{ f(x_0-r)})\sim c
\left(\frac{r}{f(x_0-r)}\right)^{N-1}$, and the same holds for
$V_A(p_0,\frac{r/2}{f(x_0+r/2)})$, we get the thesis.
 \end{proof}\medskip

 \begin{Cor}\label{2.3.4}  Let $E$ be the standard end of local
dimension $N$ and asymptotic dimension $D$ in {\rm
(\cite{Davies}}, example 16), which is homeomorphic to
$(1,\infty)\times S^{N-1}$, endowed with the metric
$ds^2=dr^2+r^{2(D-1)/(N-1)}d\o^2$, and with the volume form
$dvol=r^{D-1}drd^{N-1}\o$. Then $\ad{E}=D$.
 \end{Cor}

 Observe that $\ad{M}$ makes sense for any metric space, hence for 
 any cylindrical end, while Davies' asymptotic dimension does not. 
 Indeed let
$E := (1,\infty)\times S^1$, endowed
with the metric $ds^2=dr^2+f(r)^2d\o^2$, and with
the volume form $dvol=f(r)drd\o$, where $f(r):=\frac{d}{
dr}(r^2\log r)$. Then $\ad{E}=2$, but $vol(E_r)$ does not satisfy
one of the inequalities in (\ref{e:davies}).
 \par

 Before closing this section we observe that the notion of
standard end allows us to construct an example which shows that we
could obtain quite different results if we used
$\liminf$ instead of
$\limsup$ in the definition of the asymptotic dimension. 
 It makes use of the following function
 $$
 f(x)=
 \begin{cases}
 \sqrt{x} &  x\in[1,a_1] \\
 2+ b_{n-1}+ c_{n-1} + (x-a_{2n-1}) &
                               x\in[a_{2n-1},a_{2n}] \\
 2+ b_{n-1}+ c_{n} + \sqrt{ x-a_{2n}+1 } &
                               x\in[a_{2n},a_{2n+1}] 
 \end{cases}
 $$
 where $a_0:=0$, $a_{n}-a_{n-1}:=2^{2^n}$, $b_{n}:= \sum_{k=1}^{n} \sqrt{ 
 2^{2^{2k+1}}+1 }$, $c_{n}:= \sum_{k=1}^{n} ( 2^{2^{2k}}-1 )$,  $n\geq1$.

 \begin{Prop}\label{2.3.5} Let $M$ be the Riemannian manifold
obtained as a $\cc^\infty$ regularization of
$C\cup_\f E$, where 
$C:= \{ (x,y,z)\in\br^3 :
(x-1)^2+y^2+z^2=1,\ x\leq1 \}$, with the Euclidean metric, 
$E:=[1,\infty)\times S^1$, endowed with the metric
$ds^2=dx^2+f(x)^2d\o^2$, and with the volume form
$dvol=f(x)dxd\o$, where  $\f$ is the identification of $\{
y^2+z^2=1,\ x=1 \}$ with
$\{1\}\times S^1$. Then the volume form is a uniformly bounded
measure, $\ad{M}\geq2$ but $\underline{d}_\infty(M)\leq3/2$, where
$\underline{d}_\infty(M):=\limr\liminf_{R\to\infty}\frac{\log
n_r(B_M(x,R)) }{ \log R}$.
 \end{Prop}
 \begin{proof}
 Set $o:=(0,0,0)\in M$, then it is easy to see that, for $n\to\infty$, 
 $a_{n}\sim 2^{2^n}$, $b_{n}\sim c_{n} \sim 2^{2^{2n}}$, and 
 \begin{align*}
 area(B_M(o,a_{2n})) & \sim \frac12 a_{2n}^2 \\
 area(B_M(o,a_{2n-1}) )& \sim \frac53 a_{2n-1}^{3/2} 
 \end{align*}
 so that, calculating the limit of $\frac{\log area(B_M(o,R))}{\log R}$ 
 on the sequence $R=a_{2n}$ we get $2$, while on the
sequence $R=a_{2n-1}$ we get $3/2$. The thesis follows
easily, using Proposition \ref{1.2.14}.
 \end{proof}\medskip

 \section{A trace for open manifolds}
 \label{sec:trace}
 
 This section is devoted to the construction of a trace on (a
suitable subalgebra of) the bounded operators on $L^2(M)$, where
$M$ is an open manifold of bounded geometry. The basic idea for this construction is due to J.~Roe
 \cite{RoeOpen}, who considers regularly exhaustible open
manifolds. In our case we may (and will) restrict to exhaustions
by spheres with linearly increasing radii. Moreover we shall
regularize (three times) this trace, in order to get a
semicontinuous semifinite trace on the C$^*$-algebra of almost
local operators. As observed by Roe, this trace is strictly
related to the trace constructed by Atiyah \cite{Atiyah} in the
case of covering manifolds, and may therefore be used to define
the (0-th) Novikov-Shubin invariant for open manifolds, as we do
in Section 4.

 \subsection{The C$^*$-algebra of almost local 
 operators}\label{subsec:almostlocal}

 Recall \cite{RoeCoarse} that an operator $A\in\cb(L^2(M))$ has
finite propagation speed if there is a constant $u(A)>0$ s.t. for
any compact subset $K$ of $M$, any $\f\in L^2(M)$,  $\supp\f\subset
K$, we have $\supp A\f
\subset Pen(K,u(A)) := \{ x\in M : \d(x,K)\leq u(A) \}$. \\
 Let us denote with $\ca_0$ the set of finite propagation
speed operators. $\ca_0$ may be characterized as
follows

 \begin{Prop}\label{3.1.1} 
 \itm{i} $A\in\ca_0$ iff, for any measurable set $\O$,
$AE_\O=E_{Pen(\O,u(A))}AE_\O$, where $E_X$ is the multiplication
operator by the characteristic function of the set $X$; 
 \itm{ii} $A\in\ca_0$ iff, for any functions $\f,\ \psi\in L^2(M)$
with $\d(\supp\psi, \supp\f)\geq u(A)$, one has $(\f,A\psi)=0$.
 \end{Prop}
 \begin{proof}
 $(i)$ is obvious. \\
 $(ii)\ (\imply)$ is easy. \\
 $(ii)\ (\coimply)$ The hypothesis implies that $\supp A\psi
\subset M\setminus\supp\f$ for all $\f$  s.t.
$M\setminus\supp\f \subset Pen(\supp\psi,u(A))$. The thesis
follows.
 \end{proof}\medskip

 \begin{Prop}\label{3.1.2} The set $\ca_0$ of finite propagation
speed operators is a $^*$-al\-ge\-bra with identity.
 \end{Prop}
 \begin{proof}
 Let $K$ be a compact subset of $M$, $\f\in L^2(M)$, 
$\supp\f\subset K$, and $A,\ B\in \ca_0$. Then 
 $\supp(A+B)\f \subset \supp A\f \cup \supp B\f$, which implies
$u(A+B)\leq u(A)\vee u(B)$. 
 Moreover $\supp(AB)\f \subset Pen(\supp B\f,u(A)) \subset Pen(K,
u(A) + u(B))$, so that $u(AB) \leq u(A) + u(B)$. \\
 As $(A^*\psi,\f) = (\psi,A\f)=0$ for all $\f,\ \psi\in L^2(M)$, with $\d(\supp\psi, \supp\f)\geq u(A)$,  that is
$\supp \f \cap Pen(\supp \psi,u(A)) = \emptyset$,
 we get $\supp A^*\psi \subset
Pen(\supp\psi,u(A))$, which implies $u(A^*)\leq u(A)$, and
exchanging the roles of $A,\ A^*$, we get $u(A)= u(A^*)$.
 \end{proof}\medskip
 
 The norm closure of $\ca_0$ will be denoted by $\ca$ and will be
called the C$^*$-algebra of almost local operators. Now we show
that Gaussian decay for the kernel of a positive operator $A$ is
a sufficient condition for $A$ to belong to $\ca$. 

 \begin{Thm}\label{3.1.3} Let $M$ be an open $n$-manifold of bounded
geometry. If $A$ is a self-adjoint bounded operator on $L^2(M)$,
with kernel $a(x,y)$, and there are positive constants $c,\
\a,\ \d_0$ s.t.
 $$ 
 |a(x,y)|\leq c\ \e{-\a\d(x,y)^2},\quad \d(x,y)\geq\d_0
 $$
 then $A\in\ca$.
 \end{Thm}
 In order to prove the theorem, we need some lemmas.
 
 \begin{Lemma}\label{3.1.4} Let $A$ be a bounded self-adjoint operator
on $L^2(M)$, with measurable kernel. Then
 $$
 \|A\| \leq \sup_{x\in M} \int_M |a(x,y)|dy
 $$
 \end{Lemma}
 \begin{proof}
 Since $A$ is self-adjoint, $a(x,y)$ is symmetric, hence
 \begin{align*}
 \|A\|_{1\to1} & = \sup \{|(f,Ag)| : f\in L^\infty(M),\
\|f\|_\infty=1,\ g\in L^1(M),\ \|g\|_1=1 \} \\
 & \leq \sup_{x\in M} \int_M |a(y,x)|dy =
\|A\|_{\infty\to\infty}
 \end{align*}
 The thesis easily follows from Riesz-Thorin interpolation
theorem.
 \end{proof}\medskip

 \begin{Lemma}\label{3.1.5} Let $\f:[0,\infty)\to[0,\infty)$ be a
non-increasing measurable function. Then
 $$
 \sup_{x\in M} \int_M \f(\d(x,y))dy \leq C_n\ \int_0^\infty
\f(r)S_{c_2}(r)^{n-1}dr
 $$
 where $C_n:= \frac{n\sqrt{\pi}}{\G(n/2+1)}$, and $S_{c_2}(r):=
\frac1{\sqrt{-c_2}}\sinh(r\sqrt{-c_2})$
 \end{Lemma} 
 \begin{proof}
 From Theorem \ref{2.1.2} we get $V(x,r)\leq C_n\ \int_0^r
S_{c_2}(t)^{n-1}dt$. Then
 \begin{align*}
 \int_M \f(\d(x,y))dy & = \int_0^\infty \f(r)dV(x,r) \\
 & \leq C_n\ \int_0^\infty \f(r) S_{c_2}(r)^{n-1}dr
 \end{align*}
 where the  equality is in $e.g.$ (\cite{HeSt}, Theorem
12.46), and the inequality holds because $\f$ is non-increasing
and positive, and $V(x,0)=0$.
 \end{proof}\medskip

 \noindent {\it Proof of Theorem \ref{3.1.3}.} 
 Let $\r>\d_0$, and decompose $A=A_\r+A'_\r$, with $a_\r(x,y):=
a(x,y)\chi_{[0,\r]}(\d(x,y))$. Then $A_\r\in\ca_0$, and
$|a'_\r(x,y)|\leq c'\f(\d(x,y))$, where
 $$
 \f(r):=
 \begin{cases}
 \e{-\a\r^2} &  0\leq r<\r \\
 \e{-\a r^2} &  r\geq\r. 
 \end{cases}
 $$
 By Lemmas \ref{3.1.4}, \ref{3.1.5} we get 
 \begin{align*}
 \|A-A_\r\| = \|A'_\r\| & \leq \sup_{x\in M} \int_M
|a'_\r(x,y)|dy \\
 & \leq c'\ \sup_{x\in M} \int_M \f(\d(x,y))dy \\
 & \leq c' \int_0^\infty \f(r)S_{c_2}(r)^{n-1}dr \\
 & \leq c'\ \e{-\a\r^2} \int_0^\r S_{c_2}(r)^{n-1}dr + c''
\int_\r^\infty \e{-\a r^2+(n-1)r\sqrt{-c_2}}dr\\
 & \to 0,\quad \r\to\infty
 \end{align*}
 and the thesis follows.
 \qed\medskip
 
 Finally we conclude that $\cc_0$ functional calculus of the
Laplace operator belongs to $\ca$.

 \begin{Cor}\label{3.1.6} Let $M$ be an open manifold of bounded
geometry. Then $\f(\D)\in\ca$, for any $\f\in\cc_0([0,\infty))$.
 \end{Cor}
 \begin{proof}
 By Proposition \ref{2.1.5} and Theorem \ref{3.1.3} we obtain that
$\e{-t\D}\in\ca$, for any $t>0$. Since $\{ \e{-t\l} \}_{t>0}$
generates a dense $^*$-algebra of $\cc_0([0,\infty))$, the thesis
follows by Stone-Weierstrass theorem.
 \end{proof}\medskip

 \subsection{A functional described by J.~Roe}
 \label{subsec:roe}

 In the rest of this paper $M$ is a complete
Riemannian $n$-manifold of bounded geometry (as in Definition
\ref{2.1.1}) and of regular polynomial growth, that is
 $$
 \limr \frac{V(x,r+R) }{ V(x,r)} =1
 $$
 for all $x\in M$, $R>0$.
 \\ 
 
 Following Moore-Schochet \cite{Moore-Schochet}, we recall that an operator $T$ on $L^2(M)$ is called
locally trace class if, for any compact set $K\subset M$,
$E_KTE_K$ is trace class, where $E_K$ denotes the projection
given by the characteristic function of $K$. It is known
 that the functional $\m_T(K):=
Tr(E_KTE_K)$ extends to a Radon measure on $M$. To state next
definition we need some preliminary notions.

 \begin{Lemma}\label{3.2.1} Let $M$ be as above, and
$\m$ be a measure on $M$. Then the following are equivalent
 \itm{i} there are $x_0\in M$, $r_0>0$, $c>0$ s.t.
 $\m(B(x_0,r))\leq c\ V(x_0,r)$, $r\geq r_0$; 
 \itm{ii} $\limsup_{r\to\infty}  \frac{\m(B(x,r))}{V(x,r)}$ is
finite and independent of $x\in M$.
 \end{Lemma}
 \begin{proof}
 $(ii)\imply(i)$ is easy. \\
 $(i)\imply(ii)$ To prove the converse, we observe that, with
$\d:=\d(x_0,x)$, 
 \begin{align*}
 \frac{\m(B(x_0,r))}{V(x_0,r)} 
 & \leq \frac{\m(B(x,r+\d))}{V(x,r+\d)}\frac{V(x_0,r+2\d)}{
V(x_0,r)} \\
 & \leq \frac{\m(B(x_0,r+2\d))}{V(x_0,r+2\d)}\frac{V(x,r+3\d)}{
V(x,r+\d)} \frac{V(x_0,r+2\d)}{V(x_0,r)} 
 \end{align*}
 Taking the $\limsup_{r\to\infty}$, and making use of regular
polynomial growth, the thesis follows.
 \end{proof}\medskip

 \begin{Dfn}\label{3.2.2} A measure $\m$ on $M$ is said to be
dominated by the volume measure $vol$ in the large, denoted
$\m\prec_\infty vol$, if the equivalent conditions of the
previous Lemma hold. 
 \end{Dfn}

 Define $\cj_{0+}$ as the set of positive locally
trace class operators $T$, such that $\m_T\prec_\infty vol$, and
set $c_T:=\limsup_{r\to\infty} \frac {\m_T(B(x,r))}{V(x,r)}$.
 
 \begin{Lemma}\label{3.2.3} $\cj_{0+}$ is a hereditary (positive) cone
in $\cb(L^2(M))$.
 \end{Lemma}
 \begin{proof}
 If $T\in\cj_{0+}$, and $0\leq A\leq T$, then $Tr(BAB^*)\leq
Tr(BTB^*)$, for any $B\in\cb(L^2(M))$, and the thesis follows.
 \end{proof}\medskip
 
 Let $\o$ be a translationally invariant state on
$L^\infty([0,\infty))$, and consider the functional $\f$ on
$\cb(L^2(M))_+$ given by
 $$
 \f(A):= 
 \begin{cases}
 \o\left(\frac{\m_A(B(x,r))}{ V(x,r)}\right) &  A\in\cj_{0+} \\
 +\infty &  A\in\cb(L^2(M))_+\setminus\cj_{0+}
 \end{cases}
 $$
 Observe that the functional $\f$ is very similar to the 
functional defined by J.~Roe in \cite{RoeOpen}. Indeed regular polynomial growth
implies that $\{B(x,kr)\}_{k\in\bn}$ is a regular exhaustion
according to \cite{RoeOpen}. The further hypothesis that $\o$ is
translationally invariant will play a crucial role in our
construction.

 \begin{Prop}\label{3.2.4} 
 \itm{i} $\o$ is a generalized limit on $[0,\infty)$ 
 \itm{ii} $\f$ is a weight on $\cb(L^2(M))$ 
 \itm{iii} $\f$ does not depend on $x\in M$. 
 \end{Prop}
 \begin{proof}
 $(i)$ Since $\o(1)=1$, $\o$ vanishes on compact support
functions, therefore on $\cc_0([0,\infty))$, by continuity.
Hence, when it exists, $\lim_{t\to\infty} f(t) = \o(f)$, for
$f\in\cc_b([0,\infty))$. \\
 $(ii)$ Positivity of $\f$ is obvious, while linearity follows
from Lemma \ref{3.2.3} and the observation that $\m_{A+B} = \m_A+\m_B$.
\\
 $(iii)$ Let $x,y\in M$, $\d:=\d(x,y)$. Then, for any
$A\in\cj_{0+}$
 $$
 \m_A(B(x,r))\leq \m_A(B(y,r+\d)) \leq \m_A(B(x,r+2\d)),
 $$  
 from which it follows
 \begin{align}\label{(3.2.1)}
 \f_x(A) & := \o\left(\frac{\m_A(B(x,r)) }{ V(x,r)}\right) \notag\\
 & \leq \o\left(\frac{\m_A(B(y,r+\d)) }{
V(y,r+\d)}\frac{V(y,r+\d) }{ V(x,r)}\right) \\
 & \leq \o\left(\frac{\m_A(B(x,r+2\d)) }{
V(x,r+2\d)}\frac{V(x,r+2\d) }{ V(x,r)}\right) \notag
 \end{align}
 Since $\o$ is a generalized limit, we have
 \begin{align*}
 & \biggl| \o\left(\frac {\m_A(B(y,r+\d)) }{
V(y,r+\d)} \frac{V(y,r+\d) }{ V(x,r)} \right) -
\o\left( \frac{\m_A(B(y,r+\d)) }{ V(y,r+\d)} \right) \biggr| \\
  & \leq \limsup_{r\to\infty}\frac{\m_A(B(y,r+\d)) }{
V(y,r+\d)}\ \o\left(  \frac{V(y,r+\d) }{ V(x,r)} -1
 \right) \\
 \end{align*}
 Because of regular polynomial growth
 $$
 1\leq \lim_{r\to\infty}\frac{V(y,r+\d) }{ V(x,r)} \leq
\lim_{r\to\infty}\frac{V(x,r+2\d) }{ V(x,r)} =1 .
 $$
 Therefore, by translation invariance,
 $$
 \o\left(\frac{\m_A(B(y,r+\d)) }{ V(x,r)}\right) =
 \o\left(\frac{\m_A(B(y,r+\d)) }{ V(y,r+\d)}\right) =
 \o\left(\frac{\m_A(B(y,r)) }{ V(y,r)}\right) = \f_y(A)
 $$
 Analogously we show that
 $$
 \o\left(\frac{\m_A(B(x,r+2\d)) }{ V(x,r)}\right) =
 \o\left(\frac{\m_A(B(x,r)) }{ V(x,r)}\right) = \f_x(A)
 $$
 Then inequalities in (\ref{(3.2.1)}) read $\f_x(A)\leq \f_y(A) \leq
\f_x(A)$, and the thesis follows. 
 \end{proof}\medskip

 The algebra $\ca$, being a C$^*$-algebra, contains many unitary
operators, and is indeed generated by them. The algebra
$\ca_0$ may not, but all unitaries in $\ca$ may be approximated
by elements in $\ca_0$. Such approximants are $\d$-unitaries,
according to the following
 
 \begin{Dfn}\label{3.2.5} An operator $U\in\cb(L^2(M))$ is called
$\d$-unitary, $\d>0$, if $\|U^*U-1\|<\d$, and $\|UU^*-1\|<\d$.\\
 \end{Dfn}

  Let us denote with $\cu_\d$ the set of $\d$-unitaries in
$\ca_0$ and observe that, if $\d<1$, $\cu_\d$ consists of
invertible operators, and $U\in\cu_\d$ implies
$U^{-1}\in\cu_{\d/(1-\d)}$.

 \begin{Prop}\label{3.2.6} The weight $\f$ is $\eps$-invariant
for $\d$-unitaries in $\ca_0$, namely, for any $\eps\in(0,1)$, there is
$\d>0$ s.t., for any $U\in\cu_\d$, and $A\in\ca_+$,
 $$
 (1-\eps)\f(A) \leq \f(UAU^*) \leq (1+\eps)\f(A) .
 $$
 \end{Prop}

 \begin{Lemma}\label{3.2.7} If $T\in\cj_{0+}$, then $ATA^*\in\cj_{0+}$
for all $A\in\ca_0$.
 \end{Lemma}
 \begin{proof}
 For any $B:=B_M(x,r)$ we have
 \begin{align*}
 \frac{\m_{ATA^*}(B) }{ V(x,r)} 
 & = \frac{Tr(E_BATA^*E_B) }{ V(x,r)} \\
 & = \frac{Tr(E_BAE_{B(x,r+u(A))}TE_{B(x,r+u(A))}A^*E_B) }{ V(x,r)}
\\  
 & \leq \|A^*E_BA\| \frac{Tr(E_{B(x,r+u(A))}TE_{B(x,r+u(A))})
}{ V(x,r+u(A))} \frac{V(x,r+u(A))}{V(x,r)} \\  
 & \leq \|A\|^2 c_T\ \frac{V(x,r+u(A))}{ V(x,r)} \\  
 \end{align*}
 and, from regular polynomial growth, we get the thesis
$\m_{ATA^*}\prec_\infty vol$.
 \end{proof}\medskip
  
 \noindent {\it Proof of Proposition \ref{3.2.6}.}
 As in the proof of the previous Lemma we get
 $$
 \f(UAU^*) \leq \|U\|^2 \o\left(\frac{\m_A(B_M(x,r+u(U))) }{
V(x,r+u(U))} \frac{V(x,r+u(U))}{ V(x,r)} \right)
 $$
 which gives, as in the proof of Proposition \ref{3.2.4}$(iii)$,
 $$
 \f(UAU^*)\leq \|U\|^2 \f(A).
 $$
 Choose now $\d<\eps/2$, and $U\in\cu_\d$,
so that $U^{-1}\in\cu_{2\d}$, and $\f(UAU^*)\leq (1+\d) \f(A)< (1+\eps) \f(A)$.
 Replacing $A$ with $UAU^*$, and $U$ with $U^{-1}$, we obtain
 $$
 \f(A)\leq \|U^{-1}\|^2 \f(UAU^*)\leq (1+2\d)\f(UAU^*) < (1+\eps)\f(UAU^*)
 $$
 and the thesis easily follows.
 \qed\medskip
 
 Finally we observe that, from the proof of Lemma \ref{3.2.7} the
following is immediately obtained
 
 \begin{Lemma}\label{3.2.8} If $A\in\ca_0$ and $\|A\|\leq 1$, then
$\f(ATA^*)\leq \f(T)$, for any $T\in\ca_{0+}$. 
 \end{Lemma}

 \begin{rem}\label{3.2.9}  The use of a translation invariant state
$\o$ is the first regularization w.r.t. the original Roe's
procedure. The request of some invariance on $\o$ closely recalls
Dixmier traces versus Varga traces (see \cite{Dixmier,Varga,AGPS,GI1}), where the invariance
requirement yields a larger domain for the singular trace. With
this choice we get a much stronger property then the trace
property in \cite{RoeOpen}, namely bimodule-trace property.
Indeed, our
$\eps$-invariance for $\d$-unitaries obviously implies invariance
under conjugation with unitaries in $\ca_0$. In order to get a
trace on $\ca$, we need one more regularization, which makes
$\f|_\ca$ a semicontinuous trace on $\ca$. This procedure will be
discussed in the next subsection.
 \end{rem}

 \subsection{A construction of semicontinuous traces on C$^*$-algebras}
 \label{subsec:semicontinuous}

 In this subsection we consider a unital C$^*$-algebra $\ca$,
with a dense $^*$-subalgebra $\ca_0$. First we observe that, with
each weight on $\ca$, namely a functional
$\f_0:\ca_+\to[0,\infty]$, satisfying the property $\f_0(\l
A+B)=\l\f_0(A)+\f(B)$, $\l>0$, $A,\ B\in\ca_+$, we may associate
a (lower-)semicontinuous weight $\f$ with the following procedure
 \begin{equation}\label{e:weight}
 \f(A):= \sup \{ \psi(A): \psi\in\ca^*_+,\ \psi\leq\f_0 \}
 \end{equation}
 Indeed, it is known that \cite{Combes,Stratila}
 $$
 \f(A) \equiv \sup_{\psi\in\cf(\f_0)} \psi(A) 
 $$
 where $\cf(\f_0):= \{ \psi\in\ca^*_+ : \exists\ \eps>0,\
(1+\eps)\psi<\f_0 \}$. Moreover the following holds

 \begin{Thm}\label{3.3.1} {\rm \cite{QV}} The set $\cf(\f_0)$ is
directed, namely, for any $\psi_1,\ \psi_2\in\cf(\f_0)$, there is
$\psi\in\cf(\f_0)$, s.t. $\psi_1,\ \psi_2\leq \psi$.
 \end{Thm}

 From this theorem easily follows

 \begin{Cor}\label{3.3.2} Let $\f_0$ be a weight on the
C$^*$-algebra $\ca$, and $\f$ be defined as in (\ref{e:weight}). Then
 \itm{i} $\f$ is a semicontinuous weight on $\ca$ 
 \itm{ii} $\f=\f_0$ iff $\f_0$ is semicontinuous. \\
 The weight $\f$ will be called the semicontinuous
regularization of $\f_0$.
 \end{Cor}
 \begin{proof}
 $(i)$ From Theorem \ref{3.3.1}, $\f(A) = \sup_{\psi\in\cf(\f_0)} 
\psi(A) = \lim_{\psi\in\cf(\f_0)} \psi(A)$, whence linearity and
semicontinuity of $\f$ easily follow. \\
 $(ii)$ is a well known result by Combes \cite{Combes}.
 \end{proof}\medskip

 \begin{Prop}\label{3.3.3} Let $\t_0$ be a weight on $\ca$ which
is $\eps$-invariant by $\d$-unitaries in $\ca_0$ (as in
Proposition \ref{3.2.6}).
 Then the associated semicontinuous weight $\t$ satisfies the
same property.
 \end{Prop}
 \begin{proof}
 Fix $\eps<1$ and choose $\d\in(0,1/2)$, s.t. $U\in\cu_\d$
implies $|\t_0(UAU^*)-\t_0(A)|<\eps \t_0(A)$, $A\in\ca_+$. Then,
for any $U\in\cu_{\d/2}$, so that $U^{-1}\in\cu_\d$, and any
$\psi\in\ca^*_+$, $\psi\leq\t_0$, we get
 $$
\psi\circ adU(A)\leq
\t_0(UAU^*)\leq (1+\eps)\t_0(A),
 $$
 for $A\in\ca_+$, $i.e.$ $(1+\eps)^{-1}\psi\circ adU\leq
\t_0$.  Then 
 \begin{align*}
 \t(UAU^*) 
 & = (1+\eps) \sup_{\psi\leq\t_0} (1+\eps)^{-1}\psi\circ adU(A) \\
 & \leq (1+\eps) \sup_{\psi\leq\t_0} \psi(A) \\
 & = (1+\eps)\t(A).
 \end{align*}
 Since $U^{-1}\in\cu_\d$, replacing $U$ with $U^{-1}$ and $A$
with $UAU^*$, we get $\t(A)\leq (1+\eps)\t(UAU^*)$.
 Combining the last two inequalities, we get the result.
 \end{proof}\medskip

 \begin{Prop}\label{3.3.4} The weight $\t$ is a semicontinuous
trace on $\ca$, namely, setting $\cj_+:= \{A\in\ca_+:
\t(A)<\infty\}$, and extending $\t$ to the linear span $\cj$ of
$\cj_+$, we get
 \itm{i} $\cj$ is an ideal in $\ca$ 
 \itm{ii} $\t(AB)=\t(BA)$, for all $A\in\cj$, $B\in\ca$.
 \end{Prop}
 \begin{proof}
 $(i)$ Let us prove that $\cj_+$ is a unitary invariant face in
$\ca_+$, and it suffices to prove that $A\in\cj_+$ implies
$UAU^*\in\cj_+$, for all $U\in\cu(\ca)$, the set of unitaries in
$\ca$. Suppose on the contrary that there is $U\in\cu(\ca)$ s.t.
$\t(UAU^*)=\infty$. Then there is
$\psi\in\ca^*_+$, $\psi\leq\t_0$, s.t. $\psi(UAU^*) > 2\t(A)+2$.
Then we choose $\d<3$ s.t. $V\in\cu_\d$ implies
$\t(VAV^*)\leq 2\t(A)$, and an operator $U_0\in\ca_0$ s.t.
$\|U-U_0\|< \min\{ \frac{\d}{3}, \frac1{3\|A\|\|\psi\|} \}$. The
inequalities 
 $$
 \|U_0U_0^*-1\| = \|U^*U_0U_0^*-U^*\| \leq 
\|U^*U_0-1\|\|U_0^*\|+\|U_0^*-U^*\| \leq \d
 $$
 and analogously for $\|U_0^*U_0-1\|<\d$, show that
$U_0\in\cu_\d$.
 Then, since $|\psi(U_0AU_0^*)-\psi(UAU^*)|\leq
3\|\psi\|\|A\|\|U-U_0\|<1$, we get
 $$
 2\t(A)\geq \t(U_0AU_0^*) \geq \psi(U_0AU_0^*) \geq
\psi(UAU^*) - 1 \geq 2\t(A)+1
 $$
 which is absurd. \\
 $(ii)$ We only have to show that $\t$ is unitary invariant. Take
$A\in\cj_+$, $U\in\cu(\ca)$. For any $\eps>0$ we may find a
$\psi\in\ca^*_+$, $\psi\leq\t_0$, s.t.
$\psi(UAU^*)>\t(UAU^*)-\eps$, as, by $(i)$,
$\t(UAU^*)$ is finite. Then, arguing as in the proof of $(i)$, we
may find $U_0\in\ca_0$, so close to $U$ that
 \begin{align*}
 & |\psi(U_0AU_0^*)-\psi(UAU^*)|<\eps \\
 & (1-\eps)\t(A)\leq \t(U_0AU_0^*) \leq (1+\eps)\t(A). 
 \end{align*}
 Then
 \begin{align*}
 \t(A) 
 & \geq \frac1{1+\eps}\ \t(U_0AU_0^*) 
   \geq \frac1{1+\eps}\ \psi(U_0AU_0^*) \\
 & \geq \frac1{1+\eps}\ (\psi(UAU^*) -\eps) 
   \geq \frac1{1+\eps}\ (\t(UAU^*) -2\eps).
 \end{align*}
 By the arbitrariness of $\eps$ we get $\t(A)\geq \t(UAU^*)$.
Replacing $A$ with $UAU^*$, we get the thesis.
 \end{proof}\medskip
 
 The third regularization we need turns $\t$ into a
(lower semicontinuous) semifinite trace, namely guarantees that
 $$
 \t(A) = \sup\{ \t(B) : 0\leq B\leq A,\ B\in\cj_+ \}
 $$
 for all $A\in\ca_+$. This regularization is well known (see
$e.g.$ \cite{DixmierC}, Section 6), and amounts to represent
$\ca$ via the GNS representation $\pi$ induced by $\t$, define a
normal semifinite faithful trace $tr$ on $\pi(\ca)''$, and
finally pull it back on $\ca$, that is $tr\circ\pi$. It turns out
that $tr\circ\pi$ is (lower semicontinuous and) semifinite on
$\ca$, $tr\circ\pi\leq\t$, and $tr\circ\pi(A)=\t(A)$ for all
$A\in\cj_+$, that is $tr\circ\pi$ is a semifinite extension of
$\t$, and $tr\circ\pi=\t$ iff $\t$ is semifinite. \\
 We still denote by $\t$ its semifinite extension. As follows
from the construction, semicontinuous semifinite traces are
exactly those of the form $tr\circ\pi$, where $\pi$ is a tracial
representation, and $tr$ is a n.s.f. trace on $\pi(\ca)''$.

 \subsection{The regularized trace on the C$^*$-algebra of almost
local operators}
 \label{subsec:regularized}

 Now we apply the regularization procedure described in the
previous subsection to Roe's functional. First we observe that
$\t_0:=\f|_\ca$ is not semicontinuous.

 \begin{Prop}\label{3.4.1} The set $\cn_0:=\{T\in\ca_+:
\t_0(T)=0\}$ is not closed. In particular, there are operators
$T\in\ca_+$ s.t. $\t_0(T)=1$ but $\t(T)=0$ for any
(lower-)semicontinuous trace $\t$ dominated by $\t_0$.
 \end{Prop}
 \begin{proof}
 Recall from Theorem \ref{2.1.2}$(i)$ that there are
positive real functions $\b_1,\ \b_2$ s.t. $0<\b_1(r)\leq
V(x,r) \leq \b_2(r)$, for all $x\in M$, $r>0$, and
$\lim_{r\to0}\b_2(r)=0$. Therefore we can find a sequence
$r_n\dec0$ s.t. $\sum_{n=1}^\infty \b_2(r_n)<\infty$. Fix $o\in
M$, and set $X_n:= \{ (x_1,x_2)\in M\times M : n\leq \d(x_i,o)\leq
n+1,\ \d(x_1,x_2)\leq r_n \}$, $Y_n:= \cup_{k=1}^n X_k$,
$n\leq\infty$, and finally let $T_n$ be the integral operator whose
kernel, denoted $k_n$, is  the characteristic function of $Y_n$.
Since $k_n$ has compact support, if $n<\infty$, $\t_0(T_n)=0$. On
the contrary, since $Y_\infty$ contains the diagonal of $M\times
M$, clearly $\t_0(T_\infty)=1$. Finally
 \begin{align*}
 \|T_\infty-T_n\| 
 & \leq \sup_{x\in M} \int_M \chi_{\cup_{k=n+1}^\infty X_k}
(x,y)dy \\
 & \leq \sup_{x\in M} \sum_{k=n+1}^\infty \int_M
\chi_{X_k}(x,y)dy
\\
 & \leq \sup_{x\in M} \sum_{k=n+1}^\infty V(x,r_k) \\
 & \leq \sum_{k=n+1}^\infty \b_2(r_k) \to0. \\
 \end{align*}
 This proves both the assertions.
 \end{proof}\medskip
 
 Finally we give a sufficient criterion for a positive operator
$T$  to satisfy $\t_0(T)=\t(T)$, where $\t$ is the
semicontinuous semifinite regularization described in the
previous subsection.

 \begin{Prop}\label{3.4.2} Let $A\in\cj_{0+}$ be an integral
operator, whose kernel $a(x,y)$ is a uniformly continuous
function in a neighborhood of the diagonal in $M\times M$, namely
 \begin{equation}\label{e:uniform}
 \forall \eps>0,\ \exists \d>0\ : \d(x,y)<\d \imply
|a(x,y)-a(x,x)|<\eps 
 \end{equation}
 Then $\t_0(A)=\t(A)$.
 \end{Prop}
 \begin{proof}
 Consider first a family of integral operators $B_\d$, with
kernels 
 $$
 b_\d(x,y):= \frac{\b_1(\d)}{\b_2(\d)}\ \frac{\chi_{\D_\d}(x,y)}{
V(x,\d)},
 $$
 where $\D_\d := \{ (x,y)\in M\times M : \d(x,y)<\d \}$.
 Then $\sup_{x\in M} \int_M b_\d(x,y)dy = \frac{\b_1(\d)}{
\b_2(\d)} \leq 1$, and $\sup_{y\in M} \int_M b_\d(x,y)dx \leq
\frac{\sup_{y\in M} V(y,\d)}{\b_2(\d)} \leq 1$, which imply
$\|B_\d\|\leq1$, by Riesz-Thorin theorem. \\
 Fix $o\in M$, set $E_r$ for the multiplication operator by the
characteristic function of $B(o,r)$, and observe that
 \begin{align*}
 Tr(E_rB_\d B_\d^*E_r) 
 & = \int_{B(o,r)}dx\int_M b_\d(x,y)^2dy  \\
 & \leq \frac{\b_1(\d)}{\b_2(\d)^2}\  V(o,r) 
   \leq \frac{V(o,r)}{\b_2(\d)} \\
 \end{align*}
 Therefore $\t_0(B_\d B_\d^*)\leq \b_2(\d)^{-1}$. This implies that
$\psi_\d := \t_0(B_\d\cdot B_\d^*)$ belongs to $\ca^*_+$, and
$\psi_\d\leq \t_0$ by Lemma \ref{3.2.8}. By the results of the previous
subsection, we have $\psi_\d(A) \leq \t(A) \leq \t_0(A)$, for any
$A\in\ca_+$. \\
 Take now $A\in\ca_+$ satisfying (\ref{e:uniform}), for a pair $\eps>0$,
$\d>0$, and, setting $\b(\d):= (\frac{\b_1(\d) }{ \b_2(\d)})^2$ to
improve readability, compute
 \begin{align*}
  |Tr & (E_rB_\d AB_\d^*E_r) - Tr(E_rAE_r)| \\
 & \leq  |Tr(E_rB_\d AB_\d^*E_r) - \b(\d) Tr(E_rAE_r)| +
(1-\b(\d)) Tr(E_rAE_r) \\
 & \leq \int_{B(o,r)} dx \int_{B(x,\d)\times B(x,\d)}
b_\d(x,y)|a(y,z)-a(x,x)|b_\d(x,z) dydz \\
     & \qquad + (1-\b(\d)) Tr(E_rAE_r) \\
 & \leq  3\eps \int_{B(o,r)} dx
\int_{B(x,\d)\times B(x,\d)} b_\d(x,y)b_\d(x,z) dydz \\
     & \qquad + (1-\b(\d)) Tr(E_rAE_r) \\
 & \leq 3\eps \b(\d) V(o,r) + (1-\b(\d)) Tr(E_rAE_r) \\
 \end{align*}
 This implies $|\psi_\d(A)-\t_0(A)| \leq 3 \eps \b(\d) +
(1-\b(\d)) \t_0(A)$.
 By the arbitrariness of $\eps$ and Theorem \ref{2.1.2}$(ii)$, we get
the thesis. 
 \end{proof}\medskip

 \begin{Cor}\label{3.4.3} For any $t>0$
$\t_0(\e{-t\D}) = \t(\e{-t\D})$, where $\D$ is Laplace-Beltrami
operator.
 \end{Cor}
 \begin{proof}
 Follows from Propositions \ref{2.1.5} and \ref{3.4.2}.
 \end{proof}\medskip

 \section{Singular traces for open manifolds}
 \label{sec:singular}

\subsection{Singular traces on C$^*$-algebras}\label{subsec:singtrac}
In this subsection we shall briefly recall how to construct type II$_1$ singular 
traces on a C$^*$-algebra with a semicontinuous semifinite trace,
as is treated in \cite{GI2}. 
As it is known \cite{DixmierC}, if $\tau$ is a 
semicontinuous semifinite trace on a C$^*$-algebra $\ca$ and $\pi_\tau$ denotes the GNS
representation, there is a normal semifinite faithful trace on
$\cam:=\pi_\tau(\ca)''$ (which we still denote by $\tau$) such that
$\tau=\tau\cdot\pi_\tau$. The main problem is that while type I$_\infty$ singular
traces (like Dixmier traces, see \cite{Dixmier,AGPS}) are defined on
suitable ideals of a semifinite von~Neumann algebra $\cam$, and therefore they 
can be ``pulled back'' via $\pi_\tau$ on $\ca$, type II$_1$ singular traces, 
which are needed here, are defined on bimodules of measurable operators
affiliated to $\cam$. Then we need a notion of operator affiliated to 
$\ca$ that allows us to construct an $^*$-bimodule over $\ca$ of such 
operators and a trace on it. Moreover we need to extend $\pi_\tau$ to such a
bimodule, this extension taking values in the measurable operators affiliated to
$\cam$, and then ``pull back'' the singular traces as before. Indeed we shall see
that measurable operators affiliated to $\ca$ form what we may call a $\tau$
almost everywhere bimodule, in the sense that the usual bimodule properties will
hold only up to a zero trace projection, provided that operations are intended 
in a strong sense, as in \cite{Se}. From the technical point of view, we make use
of the ideas of Segal \cite{Se}, Nelson \cite{Ne} and Christensen \cite{Christensen}
on noncommutative integration, adapting them to the case of C$^*$-algebras. The
main problem will be the possible absence of enough projections in $\ca$, to carry
out Christensen construction, and therefore we shall construct a $^*$-algebra
containing $\ca$ on which the trace naturally extends and with enough projections
in it. In the following $(\ca,\tau)$ will be a norm closed, unital $^*$-algebra of
operators acting on a Hilbert space $\ch$ together with a semicontinuous
semifinite trace.

 \begin{Dfn} 
We say that a projection $e\in\ca''$ is {\it
essentially clopen} (w.r.t. $(\ca,\tau)$) if for all $\eps >0$, exist
$a_-,a_+\in\ca$ s.t. 
 $$
 0\leq a_-\leq e\leq a_+\leq 1,\qquad \tau(a_+-a_-)<\eps. 
 $$
 We shall denote by $\ce$ the class of $\tau$-finite essentially
clopen projections.
 \end{Dfn}
 
 \begin{Prop} {\rm \cite{GI2}} The set $\ce$ with the
operations $\vee$, $\wedge$ is a lattice.
 \end{Prop}
 
 \begin{Thm} {\rm \cite{GI2}} There exists a $^*$-subalgebra $\cc$ of $\ca''$ 
 with the following properties
 \item{$(i)$} $\ca\cup\ce\subset\cc$
 \item{$(ii)$} If $x\in\cc$, for any $\eps>0$ there exist $a_-$, 
 $a_+\in\ca$ such that
$a_-<x<a_+$ and $\tau(a_+-a_-)<\eps$
 \item{$(ii)$} The GNS representation
$\pi_\tau$ extends to a $^*$-homomorphism (still denoted by $\pi_\tau$) of $\cc$
to $\pi_\tau(\ca)''$. 
 \end{Thm}
 
 According to the preceding theorem $\cc$ is equipped
with a (positive) trace, still denoted by $\tau$, given by the pull back of the
trace on $\pi_\tau(\ca)''$, which is the unique extension of the trace on
$\ca$. The construction of $\cc$ is rather involved, indeed its elements are
not explicitly characterized, while its definition resembles that of
the enveloping Borel algebra in \cite{Ped}, therefore we shall not describe
it here. Now we pass to the definition of affiliated operators.

 \begin{Dfn} A sequence $\{e_n\}$ of essentially clopen projections is called a 
 {\it Strongly Dense Domain} (SDD) if $e_n^{\perp}$ is $\tau$-finite and
$\tau(e_n^{\perp})\to0$. We shall denote by $e$ the projection
$\sup_ne_n$.
 \end{Dfn}
 
 Let us remark that, since the trace $\tau$ is not faithful,
$e$ is not necessarily $1$. Nevertheless it is easy to show that $e^{\perp}\in\ce$
and $\tau(e)=0$. Now let us consider a linear operator $T$ acting on $\ch$. If $T$
is neither densely defined nor closed then its adjoint is a closed operator from a
proper subspace $\ck_1$ to another proper subspace $\ck_2$ of $\ch$. We shall
denote by $T^+$ the closed, densely defined operator given
by $T^+|_{\ck_1}=i_2\cdot T^*$, where $i_2$ is the embedding of $\ck_2$ into
$\ch$, and by $T^+|_{\ck_1^{\perp}}=0$. Then we denote by $T^{\nat}$ the closed
densely defined operator $(T^+)^*$. Let us recall that an operator $T$ on $\ch$ is
said to be {\it affiliated} to a von~Neumann algebra $\cam$ ($T\aff M$) if all
elements of $x\in\cam'$ send its domain into itself and $Tx\eta=xT\eta$, for any
$\eta$ in $\cd(T)$.

 \begin{Dfn} We call $\tilde\cc$ the family of
closed, densely defined operators affiliated to $\ca''$ for which there exists 
a SDD $\{e_n\}$ such that
 \itm{i} $e_n\ch\subset\cd(T)$
 \itm{ii} $e_n\ch\subset\cd(T^*)$
 \itm{iii} $Te_n\in\cc$
 \itm{iv} $T^*e_n\in\cc$. 
 \end{Dfn}
 
 If $T$, $S\in\tilde\cc$, $a\in\ca$, we consider the following 
 (strong sense) operations 
 $$
 T\oplus S:=(T+S)^{\nat},\qquad 
 a\odot T:=(a\cdot T)^{\nat},\qquad
 T\odot a:=(T\cdot a)^{\nat}. 
 $$
 We also introduce the relation of $\tau$-a.e.
equivalence, which turns out to be an equivalence relation, among operators in
$\tilde\cc$, namely $T$ is equivalent to $S$ $\t$-a.e. if there exists a common SDD $\{e_n\}$ for
$T$ and $S$ such that, setting $\ch_0:=\cup_ne_n\ch$, we have
$eT|_{\ch_0}=eS|_{\ch_0}$. We remark that, while this relation may appear
too weak, it becomes an equality as soon as the trace is faithful on $\cc$. In
fact strong sense operations too become the usual strong sense operations defined
by Segal in the case of a faithful trace on a von~Neumann algebra, therefore, as
follows by next theorem, the class of operators in $\tilde\cc$ which are $0$
a.e. are in the kernel of the extension of the GNS representation $\pi_\tau$. In
the following we shall denote by $\pi$ the GNS representation of $\ca$ associated
with the trace $\tau$, by $\cam$ the von~Neumann algebra $\pi(\ca)''$, and by
$\tilde\cam$ the algebra of measurable operators affiliated to
$\cam$.

 \begin{Thm} {\rm \cite{GI2}} The set $\tilde\cc$ is closed under strong sense
operations, and the usual properties of a $^*$-bimodule over $\ca$ hold
$\tau$-almost everywhere. Moreover the GNS representation extends to a map from
$\tilde\cc$ to $\tilde\cam$ which preserves strong sense operations.
 \end{Thm}
 
 Let us recall that, if $A\in\tilde\cam$, its distribution function and non-decreasing
rearrangement are defined as follows (cf. e.g. \cite{FK,GI1})
 \begin{align*}
 \l_A(t)&:=\t(\c_{(t,+\infty)}(|A|)) \\
 \m_A(t)&:=\inf\{s\geq0: \l_A(s)\leq t\}.
 \end{align*}
We may define the distribution function (and
therefore the associated non-decreasing rearrangement) w.r.t. $\t$ of an operator
 $A\in\tilde\cc$ as $\l_A(t)=\l_{\pi(A)}(t)$, and we get 
 $\m_A=\m_{\pi(A)}$. Then
the preimage $\ov\cc\subset\tilde\cc$ under $\pi$ of the set
$\ov\cam:=\{A\in\tilde\cam: \l_A(t_0)<\infty$ for some $t_{0}>0\}$ is an
a.e. $^*$-bimodule over $\ca$. Let us observe that, if $A\in\ov\cc$ is a positive
(unbounded) continuous functional calculus of an element in $\ca$, then
$\c_{(t,+\infty)}(A)$ belongs to $\ce$ a.e., therefore its distribution function may be
defined without using the representation $\pi$
 $$
 \l_A(t)=\t(\c_{(t,+\infty)}(A)).
 $$
 We may carry out the construction of singular traces (with respect
to $\tau$) as it has been done in \cite{GI1}. However, since only type
II$_1$ traces will be used in the following, we shall restrict to this case,
which corresponds to eccentricity at 0. 

 \begin{Dfn} An operator $T\in\ov\cc$ is
called {\it eccentric} (at 0) if either  
 $$
 \int_0^1\mu_T(t)<\infty \quad {\text{and}} \quad 
 \limsup_{t\to0}\frac{\int_0^t\m_T(s)ds}{\int_0^{2t}\m_T(s)ds}=1 
 $$
 or
 $$
 \int_0^1\mu_T(t)=\infty\quad{\text{and}} \quad
\liminf_{t\to0}\frac{\int_t^1\m_T(s)ds}{\int_{2t}^1\m_T(s)ds}=1.
 $$
 \end{Dfn}
 
 The following proposition trivially holds
 
 \begin{Prop} Let $(\ca,\t)$
be a C$^*$-algebra with a semicontinuous semifinite trace, $\p$ the associated GNS
representation, $T\in\ov\cc$, and let $X(T)$ denote the $^*$-bimodule over $\ca$
generated by $T$ in $\ov\cc$, while $X(\pi(T))$ denotes the 
$^{*}$-bimodule over $\cam$
generated by $\pi(T)$ in $\ov\cam$. Then
 \item{$(i)$} $T$ is eccentric if and only if $\pi(T)$ is
 \item{$(ii)$} $\pi(X(T))\subset X(\pi(T))$.
 \end{Prop}
 
 As in the case of von~Neumann algebras, with any eccentric operator (at 0) in 
 $\ov\cc$ we may associate a singular trace, where the word singular refers to 
 the original trace $\t$. Indeed such singular trace will vanish on bounded 
 operators. Of course singular traces may be described as the pull-back of the 
 singular traces on $\cam$ via the (extended) GNS representation. On the other 
 hand, explicit formulas may be written in terms of the non decreasing 
 rearrangement. We write these
formulas for the sake of completeness. First we observe that, by definition, if
$T\in\ov\cc$ is eccentric (at 0) there exists a pure state $\o$
on $\cc_b(0,\infty)$ which is a generalized limit in 0, namely is an extension of
the Dirac delta in 0 on $\cc[0,\infty)$ to $\cc_b(0,\infty)$, such
that
 \begin{align*}
  {\text{ if}} \quad \int_0^1\mu_T(t)<\infty & \quad{\text{then}} \quad
  \o\left(\frac{\int_0^t\m_T(s)ds}{\int_0^{2t}\m_T(s)ds}\right)=1 \\
  {\text{ if}} \quad \int_0^1\mu_T(t)=\infty & \quad {\text{then}} \quad
  \o\left(\frac{\int_t^1\m_T(s)ds}{\int_{2t}^1\m_T(s)ds}\right)=1.
 \end{align*}
 Then the singular trace associated with $\t$, $T$ and $\o$ may be written as 
 follows on the a.e. positive elements of $X(T)$, i.e. elements whose image under
 $\pi$ is positive
 \begin{align}\label{e:traces}
  \int_0^1\mu_T(t)<\infty & \quad \imply \quad 
  \t_\o(A):=\o\left(\frac{\int_0^t\m_A(s)ds}{\int_0^{t}\m_T(s)ds}\right),
    \quad A\in X(T)_+ \notag\\
  \\  
  \int_0^1\mu_T(t)=\infty & \quad \imply \quad
  \t_\o(A):=\o\left(\frac{\int_t^1\m_A(s)ds}{\int_{t}^1\m_T(s)ds}\right),
    \quad A\in X(T)_+ \notag
 \end{align}
 According to the previous analysis, some results in \cite{GI1} may 
 be rephrased as follows
 
 \begin{Thm} 
 The functionals defined in formula~{\rm(\ref{e:traces})} extend to traces
 on the a.e. $^*$-bimodule over $\ca$ $X(T)$. They also naturally extend to traces on 
 $X(T)+\ca$.
 \end{Thm} 
 
 Now, for any $T\in\ov\cc$, we set
 $$
  \a(T):=\left( \liminf_{t\to0}\frac{\log\m_T(t)}{\log\frac1{t}} \right)^{-1}
 $$
 As we shall see in the following, this number may
be considered as a generalized Novikov-Shubin invariant of $T$. A sufficient
condition for being singularly traceable (at 0) is given in terms of this
number.

 \begin{Thm}\label{Thm:singular}
 Let $T\in\ov\cc$ with $\a\equiv\a(T)$. If $\a=1$ then $T$ is
eccentric, hence singularly traceable. In general, if $\a\in(0,\infty)$ then
$T^{\a}$ is eccentric at 0.
 \end{Thm}
 \begin{proof} 
 The first statement is proved in \cite{GI3}. Then, by the properties of the 
 non-increasing rearrangement, $\m_{T^\a}(t)=\m_T(t)^\a$,  therefore
 \begin{align*}
 \a(T^\a) & = 
 \left(\liminf_{t\to0} \frac{\log\m_{T^\a}(t)}{\log\frac{1}{t}}\right)^{-1} 
  = \left(\liminf_{t\to0} \frac{\log(\m_T(t))^\a}{\log\frac{1}{t}}\right)^{-1} \\
 & = \left(\a\ \liminf_{t\to0}\frac{\log\m_T(t)}{\log\frac{1}{t}}\right)^{-1}=1
 \end{align*}
\end{proof}\medskip

 \subsection{A singular trace associated with the
Laplacian} \label{subsec:laplacian}

 In this subsection we consider an open manifold with
bounded geometry and regular polynomial growth, i.e.
the same hypotheses assumed in section 3.

 \begin{Thm}\label{4.2.1} Let $M$ be an open manifold
with bounded geometry and regular polynomial growth. Then
 $$
 \a(\D^{-1})
 =\limsup_{t\to0}\frac{\log\t(e_\D(t))}{\log t}
 =\limsup_{t\to\infty}\frac{\log \t(\e{-t\D})}{\log\frac1{t}}
 $$
 where $e_\D$ denotes the spectral family of $\D$.
 \end{Thm}
 
We need the following Lemma

 \begin{Lemma}\label{4.2.2} Let $\l:\br_+\to\br_+$ be a non-increasing,
  right continuous function, 
$\m(t):=\inf\{s\geq0:\l(s)\leq t\}$. Then
 $$
\left(\liminf_{t\to0}\frac
{\log\m(t)}{\log\frac1{t}}\right)^{-1}
=\limsup_{s\to\infty}\frac
{\log\l(s)}{\log\frac1{s}}
 $$
where the values 0 and $\infty$ are allowed.
 \end{Lemma}
 \begin{proof} 
 First recall that $\m$ is non increasing and
right continuous and that 
$\l(t)\equiv\inf\{s\geq0:\m(s)\leq t\}$. Then let $t_n\to0$ be a
sequence such that $\lim_{n\to\infty}\frac{\log\m(t_n)}{\log\frac1{t_n}}=
\liminf_{t\to0}\frac{\log\m(t)}{\log\frac1{t}}$, and let 
$t'_n:=\inf\{s\geq0:\m(s)=\m(t_n)\}
=\min\{s\geq0:\m(s)=\m(t_n)\}$ where the last equality
holds because of right continuity. Then
 $$
\liminf_{t\to0}\frac{\log\m(t)}{\log\frac1{t}}
\leq\lim_{n\to\infty}\frac{\log\m(t'_n)}{\log\frac1{t'_n}}
\leq\lim_{n\to\infty}\frac{\log\m(t_n)}{\log\frac1{t_n}}
=\liminf_{t\to0}\frac{\log\m(t)} {\log\frac{1}{t}}
 $$
namely we may replace $t_n$ with $t'_n$ to reach the
$\liminf$. Also,
$\l(\m(t'_n))=\inf\{t\geq0:\m(t)\leq\m(t'_n)\}=t'_n$,
therefore
 \begin{align*}
\liminf_{t\to0}\frac{\log\m(t)}{\log\frac{1}{ t}}
&=\lim_{n\to\infty} \frac{\log\m(t'_n)}{\log\frac{1}{ t'_n}}
=\lim_{n\to\infty} \frac{\log\m(t'_n)}{\log\frac{1}{ \l(\m(t'_n))}}\\
&\geq\liminf_{s\to0}\frac{\log s}{\log\frac{1}{ \l(s)}}
=\left(\limsup_{s\to0}\frac{\log\l(s)}{\log\frac{1}{ s}}\right)^{-1}
 \end{align*}
 For the converse inequality, let $s_n\to\infty$ be a sequence for which
$\lim_{n\to\infty}
\frac{\log\l(s_n)}{\log\frac{1}{ s_n}}=\limsup_{s\to0}
\frac{\log\l(s)}{\log\frac{1}{ s}}$. As before,
$s'_n:=\inf\{s\geq0:\l(s)=\l(s_n)\}
=\min\{s\geq0:\l(s)=\l(s_n)\}$ still brings to the
$\limsup$ and verifies $\m(\l(s'_n))=s'_n$. 
 \end{proof}\medskip

 \noindent {\it Proof of Theorem \ref{4.2.1}.} 
 First we prove the first equality
 \begin{align*}
 \a(\D^{-1})
&=\left(\liminf_{t\to0}
\frac{\log\m_{\D^{-1}}(t)}{\log\frac{1}{ t}}\right)^{-1}
=\limsup_{s\to\infty}
\frac{\log\l_{\D^{-1}}(s)}{\log\frac{1}{ s}}\\
&=\limsup_{s\to\infty}
\frac{\log\t(\c_{(s,\infty)}(\D^{-1}))}{\log\frac{1}{ s}}
=\limsup_{s\to\infty}
\frac{\log\t(\c_{(-\infty,1/s)}(\D))}{\log\frac{1}{ s}}\\
&=\limsup_{t\to0}\frac{\log\t(e_\D(t))}{\log t}
 \end{align*}
where the second equality follows by Lemma \ref{4.2.2}. For the last equality
let us set, in analogy with \cite{GS}, 
$N(\l):=\t(e_\D(\l))$, $\th(t)=\t(\e{-t\D})$. Then it
follows that $\th$ is the Laplace transform of the
Stieltjes measure defined by $N(\l)$
 $$
\th(t)=\int\e{-\l t}dN(\l),
 $$
and the last equality follows by the Tauberian theorem
contained in the appendix of \cite{GS}, provided that
we show that $\th(t)=\co(t^{-\d})$ for some $\d>0$.
 \\
On the other hand, under the assumptions of bounded
geometry, Varopoulos \cite{Varopoulos2} proved that the
heat kernel on the diagonal has a uniform
inverse-polynomial bound, more precisely, in the
strongest form due to \cite{ChavelFeldman}, we have
 $$
\sup_{x,y\in M}p_t(x,y)\leq Ct^{-1/2}
 $$
for a suitable constant $C$. Then
 $$
\th(t)=\t(\e{-t\D})=
(gen)\lim_{r\to\infty}
\frac{\int_{B(o,r)}p_t(x,x)dvol(x)}{V(o,r)}
\leq Ct^{-1/2}
 $$
which concludes the proof.
 \qed\medskip

\begin{Dfn}\label{4.2.3} Let $M$ be an open manifold
 with bounded geometry and regular polynomial growth. Then the (0-th) Novikov-Shubin
invariant of $M$ is defined as
 \begin{equation}\label{e:NS}
\a_0(M)=2\a(\D^{-1})
=2\limsup_{t\to0}\frac{\log(N(t))}{\log t}
=2\limsup_{t\to\infty}
\frac{\log(\th(t))}{\log\frac{1}{ t}}.
 \end{equation}
 \end{Dfn}
 
 \begin{rem}
 We have chosen J.~Lott's normalization \cite{Lott} for the Novikov-Shubin 
 number $\a_{0}(M)$ because Laplace operator is a second order 
 differential operator, and this normalization gives the equality 
 between $\a_{0}(M)$ and the asymptotic dimension of $M$, cf. 
 Theorem \ref{4.3.2}. \\
 Our choice of the $\limsup$ in (\ref{e:NS}), in contrast with 
 J.~Lott's choice \cite{Lott}, is motivated by our interpretation 
 of $\a_{0}(M)$ as a dimension. On the one hand, it is 
 compatible with the classical properties of a dimension as stated in 
 Theorem \ref{Thm:adim}, cf. also Remark \ref{1.2.15}, on the other 
 hand, a noncommutative measure corresponds to such a dimension via a 
 singular trace, according to Theorem \ref{Thm:singular}, cf. 
 \cite{GI3}.  
 \end{rem}
 
\begin{Cor}\label{4.2.4} Let $M$ be an open manifold with
bounded geometry and regular polynomial
growth. Then there exists a singular trace on 
$\ov\cc$ which is finite on the $^*$-bimodule over $\ca$
generated by $\D^{-\a_0(M)/2}$.
 \end{Cor}

 \begin{rem}\label{4.2.5} This singular trace is the global
(or asymptotic) counterpart of the Wodzicki residue, in
the form of Connes, namely it is a singular trace which is
finite exactly on the operators with a prescribed
asymptotic behavior. Such an asymptotic behavior is
that of a suitable power of the Laplace operator,
i.e. that of a geometric pseudo-differential operator
with a suitable order. The problem is that such an order
seems to depend on the trace on $M$, which in turn
depends on a dilation invariant limit procedure.
Moreover, in the case of the local singular trace, such
an order has a geometric meaning, is indeed the
dimension of the manifold. These two questions will be
completely solved in the next subsection, though
under more stringent hypotheses on the manifold. 
 \end{rem}

 \subsection{The asymptotic dimension and the 0-th
Novikov Shubin invariant} \label{subsec:NSinvariant}

In this subsection, besides bounded geometry and regular polynomial 
growth, we shall also 
assume the {\it isoperimetric
inequality} which was the subject of Theorem \ref{2.1.7}, namely
there are $\a$, $\b>0$ s.t. for all $x\in M$, $r>0$, and
all regions $U\subset B(x,r)$, the first Dirichlet
eigenvalue of $-\D$ in $U$, $\l_1(U)$, satisfies
 $$
 \l_1(U) \geq \frac{\a }{ r^2} \left( \frac{V(x,r)}{ vol(U)}
\right)^\b .
 $$

 First we observe that in this case the volume of the
balls of a given radius	is uniformly bounded.

 \begin{Lemma}\label{4.3.1}	If the previous	hypotheses hold, then
 $$
  \g^{-1}\leq\frac{V(x,r)}{V(y,r)}\leq\g.
 $$
 \end{Lemma}
 \begin{proof} The Lemma easily	follows	by a result	of
Grigor'yan (\cite{Grigoryan94}, Proposition 5.2),	where it is	shown that
the	isoperimetric inequality above implies the existence
of a constant $\g$ such	that
 $$
\g^{-1}\left(\frac{R}{r}\right)^{\a_{1}} \leq \frac{V(x,R)}{V(y,r)}
\leq \g\left(\frac{R}{r}\right)^{\a_{2}}
 $$
 for some positive constants $\a_{1},\ \a_{2}$, for any $R\geq r$, and 
 $B(x,R)\cap B(y,r)\neq\emptyset$.
 \end{proof}\medskip

 \begin{Thm}\label{4.3.2} Let $M$ be as above. Then the asymptotic dimension of $M$ coincides
with the 0-th Novikov-Shubin invariant, namely
$d_\infty(M)=\a_0(M)$.
 \end{Thm}
 \begin{proof} First, from Theorem \ref{2.1.7} and the previous Lemma, we get
 \begin{align*}
 \frac{C\g^{-1}}{V(o,\sqrt{t})}
&\leq \frac {\int_{B(o,r)} \frac{C}{V(x,\sqrt{t})} dvol(x)} {V(o,r)}
\leq \frac{\int_{B(o,r)}p_t(x,x)dvol(x)}{V(o,r)}\\
&\leq \frac{\int_{B(o,r)} \frac{C'}{V(x,\sqrt{t})} dvol(x)} {V(o,r)}
\leq  \frac{C'\g}{V(o,\sqrt{t})}
 \end{align*}
therefore, by definition of the trace $\t$,
 $$
 \frac{C\g^{-1}}{V(o,\sqrt{t})}
 \leq \t(\e{-t\D})
 \leq \frac{C'\g}{V(o,\sqrt{t})}
 $$
hence, finally,
 \begin{align*}
 d_\infty(M) & = 2\limsup_{t\to\infty}\frac{\log(V(o,t))}{2\log t}
 = 2\limsup_{t\to\infty}\frac{\log(C\g^{-1} V(o,\sqrt{t})^{-1})}{\log\frac{1}{t}}\\
 \leq \a_0(M) & = 2\limsup_{t\to\infty}\frac
{\log\t(\e{-t\D})}{\log\frac{1}{t}}
\leq2\limsup_{t\to\infty}\frac
{\log(C'\g V(o,\sqrt{t})^{-1})}{\log\frac{1}{t}}\\
 & = 2\limsup_{t\to\infty}\frac
{\log(V(o,t))}{2\log t} = d_\infty(M)
 \end{align*}
\end{proof}\medskip

 \begin{rem}\label{4.3.3} On the one hand the previous result
shows that the 0-th Novikov-Shubin invariant is
intrinsically defined, since it coincides with a
rough-isometry invariant. On the other hand, the singular
trace described in the previous subsection is finite (and
non trivial) on the geometric operator $\D^{-\a/2}$,
namely on a pseudo-differential operator of degree
$-d_\infty$. In this case too such a degree plays the
role of a dimension, and more precisely coincides with
the asymptotic dimension of the manifold.
 \end{rem}
 
 \begin{ack}
 We would like to thank M.~Shubin and P.~Piazza for conversations. We 
 also thank I.~Chavel, E.B.~Davies, A.~Grigor'yan, P.~Li, L.~Saloff-Coste 
 for having suggested useful references on heat kernel estimates. 
 \end{ack}


\end{document}